\documentclass{amsart}
\usepackage{amsmath,amsfonts,amsthm,amssymb,
amsthm,graphicx,mathtools,amscd}
\usepackage{hyperref}

\usepackage[mathscr]{eucal}
\usepackage{graphics}
\usepackage{color}
\usepackage{epsfig}

\usepackage[nocompress]{cite}

\DeclareSymbolFont{iwonaletters}{OML}{iwona}{m}{it}
\DeclareMathSymbol{\bdel}{\mathalpha}{iwonaletters}{"E}

\DeclareMathOperator*{\esssup}{ess\,sup}

\theoremstyle{plain}
\newtheorem{theorem}{Theorem}[section]
\newtheorem{lemma}[theorem]{Lemma}

\newtheorem{definition}[theorem]{Definition}

\newtheorem{proposition}[theorem]{Proposition}

\newtheorem{assumption}[theorem]{Assumption}

\newtheorem{remark}[theorem]{Remark}

\newcommand{\beginsec}{
\setcounter{equation}{0}
}

\newcommand{\eps}{\varepsilon}
\newcommand{\ph}{\varphi}

\newcommand{\al}{\alpha}

\newcommand{\gam}{\gamma}
\newcommand{\kap}{\kappa}
\newcommand{\sig}{\sigma}
\newcommand{\del}{\delta}
\newcommand{\om}{\omega}

\newcommand{\Del}{\mathnormal{\Delta}}
\newcommand{\Th}{\mathnormal{\Theta}}
\newcommand{\La}{\mathnormal{\Lambda}}

\newcommand{\Sig}{\mathnormal{\Sigma}}

\newcommand{\Om}{\mathnormal{\Omega}}

\newcommand{\N}{{\mathbb N}}
\newcommand{\Q}{{\mathbb Q}}
\newcommand{\R}{{\mathbb R}}

\newcommand{\Z}{{\mathbb Z}}

\newcommand{\EE}{{\mathbb E}}
\newcommand{\E}{{\mathbb E}}

\newcommand{\PP}{{\mathbb P}}

\newcommand{\calB}{{\mathcal B}}

\newcommand{\calF}{{\mathcal F}}

\newcommand{\calI}{{\mathcal I}}

\newcommand{\calL}{{\mathcal L}}
\newcommand{\calM}{{\mathcal M}}

\newcommand{\calQ}{{\mathcal Q}}

\newcommand{\calR}{{\mathcal R}}
\newcommand{\calS}{{\mathcal S}}
\newcommand{\calT}{{\mathcal T}}

\newcommand{\frX}{\mathfrak{X}}
\newcommand{\fra}{\mathfrak{a}}
\newcommand{\frb}{\mathfrak{b}}
\newcommand{\frc}{\mathfrak{c}}

\newcommand{\frp}{\mathfrak{p}}
\newcommand{\frs}{\mathfrak{s}}

\newcommand{\lan}{\langle}
\newcommand{\ran}{\rangle}

\newcommand{\oo}{\overline}
\newcommand{\supp}{{\rm supp}}

\newcommand{\id}{{\tt id}}
\newcommand{\w}{\wedge}

\newcommand{\pl}{\partial}

\newcommand{\To}{\Rightarrow}

\newcommand{\iy}{\infty}
\newcommand{\up}{\uparrow}
\newcommand{\loc}{{\rm loc}}
\newcommand{\cadlag}{c\`adl\`ag }

\newcommand{\noi}{\noindent}
\newcommand{\ds}{\displaystyle}
\newcommand{\LE}{\preccurlyeq}
\newcommand{\DEL}{{\hat\delta}}
\newcommand{\vhat}{\widehat}

\newcommand{\MOD}{ {\rm\ mod\ }}

\newcommand{\INJ}{{\rm INJ}}
\newcommand{\REM}{{\rm REM}}

\begin{document}

\title[Free boundary problem and particle systems]{
A weak formulation of free boundary problems and its application to
hydrodynamic limits of particle systems with selection
}

\author{Rami Atar}
\address{Viterbi Faculty of Electrical and Computer Engineering
\\
Technion, Haifa, Israel
} 



\begin{abstract}

A weak formulation for a class of parabolic free boundary problems
(FBP) is proposed
that does not involve the notion of a free boundary
but reduces to a FBP when classical solutions exist.
It is aimed at hydrodynamic limits (HDL) of particle systems
with selection in circumstances where the macroscopic model
does not possess (or is hard to prove to possess) a regular
free boundary in the classical sense.
The formulation involves the macroscopic density of particles
and a measure that accounts for selection.
It consists of a second order parabolic equation satisfied by the density
and driven by the measure, coupled with a complementarity condition satisfied
by the density-measure pair. The approach is applied to
an injection-branching-selection particle system of diffusion on $\R$
under arbitrarily varying injection and removal rates, for which the corresponding FBP
is not in general known to be classically solvable.
The HDL is characterized as the unique solution to the weak formulation.
The proof of convergence is based on PDE uniqueness,
which in turn relies on the barrier method.
\end{abstract}

\maketitle

\section{Introduction}

\subsection{Background and motivation}

In particle systems with spatial selection,
particles living in $\R$ undergo motion, branching or injection, and
selection. The last term refers to keeping the population size
constant by removing, upon appearance of a new particle,
the particle whose position is smallest among all particles
(henceforth, the leftmost particle).
The first such systems were proposed in \cite{bru06, bru07}
as models for natural selection in the evolution of a population,
where the position of a particle
represents the degree of fitness of an individual
to its environment.
A series of papers culminating in the monograph \cite{de-masi-book}
studied a variety of related models motivated by particles interacting topologically
(since, macroscopically,
removals occur at the boundary of the configuration)
and by particle systems in contact with current reservoirs.
At the hydrodynamic limit (HDL) these models give rise to
free boundary problems (FBP).
Rigorously establishing the HDL--FBP relation requires control over
regularity of the free boundary (such as $C^1$ or sometimes $C$).
We are motivated by questions of characterizing HDL
by PDE in cases where existing techniques might fall short of
yielding free boundary regularity.
This may happen, in particular, when the constant population size
assumption is dropped and injection and removal rates
vary at the macroscopic scale. The first goal of this paper is to introduce
a weak formulation of FBP that does not involve
the notion of a free boundary, and at the same time reduces to a classical
FBP when classical solutions exist.

To put these questions in context,
consider the {\it $N$-particle branching Brownian motion} ($N$-BBM) in dimension 1,
first studied in \cite{mai16}, which consists of $N$ particles performing Brownian motion (BM)
independently, each branching into two at rate $1$.
When branching occurs, the leftmost particle in the configuration is removed.
The initial particle positions are drawn independently according to a probability
measure $\xi_0$. Let $\xi^N_t=\sum_{i}\bdel_{X_i(t)}$ denote the configuration measure at time $t$,
where $X_i(t)$ are the locations of the $N$ particles alive at time $t$.
Here and throughout, $\bdel_x$ denotes the Dirac measure at $x$.
Throughout this paper, the bar notation will stand for normalization by $N$;
in particular, $\bar\xi^N_t=N^{-1}\xi^N_t$.
The corresponding FBP is to find a pair $(u,\ell)$, $\ell\in C((0,\iy):\R)$,
$u\in C(\R\times(0,\iy):[0,\iy))\cap C^{2,1}(\{(x,t):t>0,x>\ell_t\})$ such that
\begin{equation}\label{401}
\begin{cases}
\pl_tu-\frac{1}{2}\pl_x^2u=u & x>\ell_t,
\\
u=0 & x\le\ell_t,\\
u(x,t)dx\to\xi_0(dx) & \text{weakly as } t\downarrow0,\\
\int_\R u(x,t)dx=1.
\end{cases}
\end{equation}
It was shown in \cite{de-masi-nbbm} (for $\xi_0$ possessing a density)
that the process $\bar\xi^N_t$ has a deterministic limit, characterized
in terms of barriers (see \S \ref{sec3}).
Under the assumption that \eqref{401} has a classical
solution and that the free boundary $\ell$ is $C^1$,
it was further shown that
the limit process has a density given by the unique solution to \eqref{401}.
In \cite{ber19} it was then shown, for general $\xi_0$,
that \eqref{401} has a unique classical solution
and that the limit of $\bar\xi^N$ has a density given by $u$
(with $\ell$ only in $C$).
In \cite{de-masi-book}, a model we will refer to as
the {\it injection-selection} model was studied,
in which a collection of $N$ Brownian particles living in $\R_+$ and reflecting
at the origin, is subject to injection of new particles at the origin,
at times determined by a rate-$c_0N$ exponential clock, $c_0>0$ a constant.
Upon each injection, the rightmost particle is removed (i.e., the one whose location
on $\R$ is greatest).
The corresponding FBP is to find $(u,\ell)$ such that
\begin{equation}\label{402}
\begin{cases}
\pl_tu-\frac{1}{2}\pl^2_xu=0 & 0<x<\ell_t,
\\
u=0 & x\ge\ell_t,\\
-\frac{1}{2}\pl_xu(0,t)=-\frac{1}{2}\pl_xu(\ell_t,t)=c_0,\\
u(\cdot,0)=u_0,
\end{cases}
\end{equation}
with $u_0$ an initial density.
It is shown there that the HDL exists and possesses a density.
Moreover, to overcome questions of regularity of
the free boundary, a weak formulation of solutions to \eqref{402} is proposed there,
defined via approximations by local classical solutions,
and it is proved that such a solution uniquely exists and is equal to the
aforementioned HDL density (see \S\ref{sec16} for more details).

The context in which the weak formulation is presented, in \S\ref{sec12} below, is
an {\it injection-branching-selection} system of diffusion processes on $\R$,
which extends both the aforementioned ones (except the minor detail that,
in \cite{de-masi-book}, particles live in $\R_+$). In this model,
the mass conservation condition is abandoned, and the rates of injection
and removal of mass may vary. Such a perturbation has dramatic consequences on
the macroscopic model to the extent that they may lead to high degree of
free boundary irregularity (see Remark~\ref{rem00}).
It is not clear whether the current toolboxes of either the classical
solution approach of \cite{ber19} or the weak solution approach of
\cite{de-masi-book} can potentially cover such scenarios.
We will show that the weak formulation introduced here does.

Our second goal has to do with the applicability of the PDE uniqueness approach
to studying HDL, which consists of showing that all limit laws are supported
on solutions to a PDE that possesses a unique solution.
The use of this approach has been missing from the literature on the subject,
precisely due to difficulties regarding free boundary regularity;
we refer to \cite{de-masi-excl} for a discussion of this point.
We will show that the weak formulation fills this gap
at least insofar as the injection-branching-selection model is concerned.

\subsection{Injection-branching-selection and weak formulation}\label{sec12}

A brief description of the model is as follows.
Brownian particles on the line, whose initial
number is $N$, branch at rate $\kap\ge0$. In addition,
injections occur according to a given point process,
and removals occur at the left edge, with
their number up to time $t$ denoted by $J^N_t$.
The space-time injection and removal locations are
denoted by $\INJ^N_i$ and $\REM_i^N$, $i\in\N$, respectively,
and are encoded in random measures on $\R\times\R_+$, namely
\begin{equation}\label{403}
\al^N(dx,dt)=\sum_{i\in\N}\bdel_{{\rm INJ}^N_i}(dx,dt),
\qquad
\beta^N(dx,dt)=\sum_{i\in\N}\bdel_{{\rm REM}^N_i}(dx,dt).
\end{equation}
Our scaling assumption is that
$(\bar\xi^N_0,\bar\al^N,\bar J^N)\to(\xi_0,\al,J)$ in probability,
where the latter is a deterministic tuple, and $J$ is absolutely continuous,
nondecreasing and null at zero
(recall that the bar notation stands for normalization by $N$).

We can now provide a formal derivation of a PDE formulation
that does not involve a free boundary. The fact that particles are always
removed from the left side of the configuration can be expressed
as a condition on $(\bar\xi^N,\bar\beta^N)$, namely
\begin{equation}\label{404}
\bar\beta^N(\{(x,t)\in\R\times\R_+:\bar\xi^N_t(-\iy,x)>0\})=0.
\end{equation}
A key point is that $\bar\beta^N$ plays an additional role in the model, namely
it drives the dynamics. Let us assume that
in some sense $(\bar\xi^N,\bar\beta^N)\to(\xi,\beta)$
as $N\to\iy$,
and moreover, that $\xi_t$ has a density $u(\cdot,t)$
for each $t$. Then the macroscopic
dynamics should satisfy
\[
\pl_tu-\frac{1}{2}\pl_x^2u-\kap u=\al-\beta.
\]
Thus one is led to the following problem formulation.
Let data $(\xi_0, \al, J)$ be given. Denote $I_t=\al(\R\times[0,t])$.
Assume that the macroscopic total mass remains positive,
namely that if $m_t=1+\kap \int_0^t m_sds +I_t-J_t$ then
$m_t>0$ for all $t$.
Find $(u,\beta)$, $u$ nonnegative, such that
\begin{equation}\label{405}
\begin{cases}
\ds
\pl_tu-\frac{1}{2}\pl_x^2u-\kap u=\al-\beta
\\
\ds
\beta(U>0)=0 \qquad \text{where} \qquad U(x,t)=\int_{-\iy}^xu(y,t)dy
\\
\beta(\R\times[0,t])=J_t
\\
u(\cdot,0)=\xi_0.
\end{cases}
\end{equation}
The precise definition of solutions to a second order parabolic equation
with measure-valued r.h.s.\ and measure initial condition
is given in \S\ref{sec2}.
We will refer to this as the {\it weak FBP formulation},
and to the condition $\beta(U>0)=0$ as the {\it complementarity condition}.

Our main result, Theorem \ref{th2}, states that, under mild assumptions
on $\al$ and $J$, there exists a unique solution $(u,\beta)$ to \eqref{405},
and, moreover, $(\bar\xi^N,\bar\beta^N)\to(\xi,\beta)$ in probability, where
$\xi_t(dx)=u(x,t)dx$. The result is stated in a broader set up in which the particles
follow a diffusion process on the line.
To recapitulate, this formulation circumvents the non-trivial obstacle
of determining conditions for existence of a free boundary as a 
regular trajectory
and related convergence issues, and, moreover,
makes it possible to argue via PDE uniqueness.

\subsection{Related work}\label{sec16}

\ 

\vspace{.5em}

\noi{\it Particle systems with selection and related models}.
A model for motionless non-locally branching particles
with selection was studied in \cite{dur11}, and its HDL was proved
to be given in terms of an integro-differential FBP.
HDL for a model that involves injection and selection
was studied in \cite{de-masi-top}, where particles perform random walks on
$[0,N]\cap\Z$.
A variant of the $N$-BBM, in which branching is nonlocal,
was studied in \cite{de-masi-nbbm2},
where the HDL was proved to exist
with explicit bounds on the rates.
The characterization of the limit as the solution of a FBP was also proved
conditionally on existence of a classical solution to the latter, but
existence is not known in general.
$N$-BBM in higher dimension with a radially symmetric
fitness function was studied in \cite{ber21, ber22bee}.
Recently, in \cite{klu23}, the HDL of a system of Brownian particles with
selection was characterized via the inverse first-passage time problem.

There are both formal and rigorous relations between FBP \eqref{401} and
\eqref{402} and the Stefan FBP \cite{de-masi-book, ber21}.
The latter was obtained as
limits of variants of the symmetric simple exclusion process (SSEP)
in \cite{landim98, landim06}, as well as the limit of interacting diffusions
with rank-dependent drift in \cite{dem19}.
In \cite{de-masi-excl}, a SSEP
with birth of the leftmost hole and death of the rightmost particle was considered,
and convergence at the hydrodynamic scale was proved;
the question of a rigorous connection to a FBP was left open.

In addition to proving HDL, some of the aforementioned papers have
studied the long time behavior of the macroscopic dynamics
and the interchange of the $N$ and the $t$ limits
\cite{de-masi-nbbm, de-masi-book, de-masi-excl,
dur11, ber21, ber22bee}.

\vspace{.5em}

\noi{\it On earlier weak formulations}.
Relaxed solutions to \eqref{402} were proposed in \cite{de-masi-book},
defined by the limit of a sequence of classical solutions to
FBP with perturbed data, as the perturbation size tends to zero.
The existence of these classical solutions was
proved via local existence to the Stefan problem
with piecewise $C^1$ free boundary.
To the best of our knowledge, this idea has been implemented
only for the model studied in \cite{de-masi-book},
where injection and removal rates are constant.
The paper \cite{delarue22} introduced a probabilistic reformulation of the
Stefan FBP and used it to define solutions beyond singularities,
known to occur in the supercooled case of the problem. While the equations are related,
this formulation does not directly apply to the FBP considered here.

\vspace{.5em}

\noi{\it On the barrier method.}
The use of barriers was introduced into the subject
in \cite{de-masi-top} and \cite{de-masi-excl} (for particle systems
with topological interaction, and, respectively, SSEP
with free boundaries). In this context, barriers are discrete versions
of the macroscopic dynamics defined by a Trotter procedure,
that bound below and above solutions to a FBP.
They have been used to obtain uniqueness by
showing that there can be at most one element separating all
lower barriers from all upper barriers. Similar ideas have been used at the level of
stochastic particle systems to provide a.s.\ bounds on particle system dynamics.
The use of deterministic barriers to proving uniqueness
of a relaxed FBP solution first appeared in \cite{de-masi-book}.
Variants of the method have since been used in several papers including
\cite{de-masi-nbbm, de-masi-nbbm2, ber19}.

A key to proving that barriers form bounds on FBP
solutions is a Feynman-Kac representation.
This tool is missing in the generality at which \eqref{405} is considered in this paper,
for reasons having to do with the irregularity of the free boundary.
Consequently, the structure of, and argument behind,
the (lower) barriers developed here differs from that in the above references.
This issue is discussed in Remark \ref{rem00} and \S\ref{sec30++}.

\subsection{Organization of the paper}

In \S \ref{sec2}, the injection-branching-selection model is constructed,
the weak FBP formulation is defined, and the main result, Theorem \ref{th2}, is stated.
The remaining sections provide the proofs.
In \S\ref{sec3}, the barrier method is used to prove PDE uniqueness,
starting from properties of mild solutions in \S\ref{sec30},
and properties of operators required for the construction of the barriers
in \S\ref{sec30+}. In \S\ref{sec30++}
the construction of barriers in earlier work is described,
and a sketch of their use in this paper is provided.
The barriers are constructed in \S\ref{sec31} and \S\ref{sec32}.
The proof of uniqueness is completed in \S\ref{sec33}.
The convergence is proved in \S\ref{sec4}, where \S\ref{sec41}
studies the measure-valued prelimit process
(Lemma \ref{lem4}), establishes tightness of $(\bar\xi^N,\bar\beta^N)$
(Lemma \ref{lem6}) and measurability of limit densities
(Lemma \ref{lem5}). \S\ref{sec42}
shows that the complementarity condition is preserved
under the limit (Lemma \ref{lem7}), and finally, \S\ref{sec43} completes the proof of
the main result.

\subsection{Notation}

Denote by $\N$ (resp., $\Z_+$) the set of positive (resp., nonnegative) integers.
For $N\in\N$, $[N]:=\{1,2,\ldots,N\}$.
Let $\mathbb{B}_r(x)=\{y\in\R:\|y-x\|\le r\}$.
Denote by $\calM(\R)$ the space of finite signed Borel measures on $\R$
endowed with the topology of weak convergence.
Let $\calM_1(\R)\subset\calM_+(\R)\subset\calM(\R)$
denote the subsets of probability and, respectively positive measures,
and give them the inherited topologies. Denote $\R_+=[0,\iy)$ and
let $\calM_{\rm loc}(\R\times\R_+)$ be the
space of signed Borel measures on $\R\times\R_+$
that are finite on $\R\times[0,T]$
for every $T$ and give it the topology of weak convergence
on $\R\times[0,T]$ for every $T$. Similarly, let
$\calM_{+,\rm loc}(\R\times\R_+)\subset\calM_{\rm loc}(\R\times\R_+)$
be the subspace of positive measures with inherited topology.
For $X=\R$ or $\R\times\R_+$, for $\mu\in\calM_+(X)$,
denote the total mass by $|\mu|=\mu(X)$ and the support
by $\supp\,\mu$. For $\mu,\nu\in\calM_+(X)$,
write $\mu\sqsubset\nu$ if $\mu(A)\le\nu(A)$ for all measurable
$A\subset X$.

For $u,v\in\calB(\R,\R)$ (Borel measurable) and $\xi\in\calM_+(\R)$,
denote $\lan u,\xi\ran=\int u\, d\xi$ and $(u,v)=\int uv\,dx$.
For $\xi\in\calM_+(\R)$, $u\in L_1(\R)$ and an interval, say $[a,b]$, use
$\xi[a,b]$ and $u[a,b]$ as shorthand for $\xi([a,b])$, and, respectively,
$\int_a^b u\,dx$.

For $p\in[1,\iy]$, abbreviate $L_p(\R)$ to $L_p$, and
for $u\in L_p$ denote $\|u\|_p= \|u\|_{L_p}$. For $p,q\in[0,\iy]$,
let $L_{p,\rm loc}(\R_+,L_q)$ denote the linear space of
functions from $\R_+$ to $L_q$ that are $p$-integrable on $[0,T]$
for every $T$, i.e., $\int_0^T\|u(\cdot,t)\|_q^pdt<\iy$ if $p\in[1,\iy)$ and
$\esssup\{\|u(\cdot,t)\|_q:t\in[0,T]\}<\iy$ if $p=\iy$, equipped with the corresponding norm.
Define $L_{p,\rm loc}((0,\iy),L_q)$ similarly, with $[0,T]$ replaced by $[T_1,T_2]$,
$0<T_1<T_2<\iy$.
A member $u=u(x,t)$ of $L_{p,\rm loc}(\R_+,L_q)$ is said to be a.e.\ non-negative if,
for a.e.\ $t$, it is non-negative for a.e.\ $x$.

For $(X,d_X)$ a Polish space let $C(\R_+,X)$ and $D(\R_+,X)$
denote the space of continuous and, respectively, \cadlag paths,
endowed with the topology of uniform convergence on compacts
and, respectively, the Skorohod $J_1$ topology.
Let $C^\up(\R_+,\R_+)$ denote the subset of $C(\R_+,\R_+)$
of nondecreasing functions that vanish at zero.
For $I\in C^\up(\R_+,\R_+)$, denote by $dI_t$
the corresponding Stieltjes measure on $\R_+$.
Denote by $AC^\up(\R_+,\R_+)$ the subset of $C^\up(\R_+,\R_+)$
of absolutely continuous functions.
For $\rho\in(0,1]$ let $C^\rho(\R_+,\R)$ denote the space of
$\rho$-H\"older continuous functions starting at zero.
Denote by $C_c^\iy(X)$ the space
of compactly supported smooth functions on $X$ when $X=\R$ or $\R\times\R_+$.
For $f:\R_+\to X$ denote
\[
w_{[T_1,T_2]}(f,\del)=\sup\{d_X(f(s),f(t)):T_1\le s\le t\le (s+\del)\w T_2\},
\]
and $w_{T}=w_{[0,T]}$.
For $(Y,|\cdot|)$ a normed space and $f:\R_+\to Y$, denote
\[
\|f\|^*_{[T_1,T_2]}=\sup\{|f(s)|:t\in[T_1,T_2]\}
 \quad\text{and}\quad
\|f\|^*_T=\|f\|^*_{[0,T]}.
\]

The term {\it with high probability} (w.h.p.)
means `away from an $N$-dependent event whose probability
tends to zero as $N\to\iy$'.
The symbol $c$ denotes a positive constant
whose value may change from one expression to another.

\section{Particle system model, weak formulation, and main result}\label{sec2}
\beginsec

\subsection{Particle system construction}\label{sec211}

First we describe the motion that individual particles perform,
namely a diffusion process with coefficients $\frb$ and $\frc$
satisfying
\begin{assumption}\label{assn0}
One has $\frb\in C^1(\R)$ and $\frc\in C^2(\R)$ with $\frb$,
$\frc$ and its derivative $\frc'$ bounded
and $\frc$ bounded away from zero.
\end{assumption}
Given a one-dimensional BM $B$, denote by $\frX(x,s,B)$
the unique strong solution $\{X_t:t\in[s,\iy)\}$ to the SDE
\begin{equation}\label{407}
X_t=x+\int_s^t\frb(X_\theta)d\theta+\int_s^t\frc(X_\theta)dB_\theta,
 \qquad t\in[s,\iy).
\end{equation}

The particle system, defined on a probability space $(\Om,\calF,\PP)$,
is indexed by $N$, where $N$ is the initial number of particles.
The particles are indexed by the set $\calS=\N\times\Z_+$, where
particles $i=(j,0)$ are roots of family trees, and particles
$i=(j,k)$, $k\ge1$ are descendants of $(j,0)$.
Below, items (S.1)--(S.5) list stochastic primitives of the model,
and (S.6) states a condition they satisfy.
\begin{enumerate}
\item[S.1]
A collection $\{x^i(N): i=(j,0),\ j\in[N]\}$ of real-valued random variables
representing the initial positions of the particles in the initial configuration.
For such $i$ set $\sig^i(N)=0$, expressing the fact that these
particles are present in the system at time $0$.
\item[S.2]
A collection $\{(x^i(N),\sig^i(N)): i=(j,0),\ j=N+1,N+2,\ldots\}$ of
$\R\times(0,\iy)$-valued random variables representing
the initial space-time positions of injected particles,
ordered by injection time, assumed distinct:
$\sig^{(N+1,0)}(N)<\sig^{(N+2,0)}(N)<\cdots$.
\item[S.3]
A collection $\{B^i:i\in\calS\}$ of mutually independent BM,
driving the motion of the corresponding particles.
\item[S.4]
A collection $\{\pi^i:i\in\calS\}$ of mutually independent
rate-$\kap$ Poisson processes, where $\kap\ge0$ is the branching rate,
determining the times a living particle gives birth.
\item[S.5]
A sequence $0<\eta^1(N)<\eta^2(N)<\cdots$ of removal attemp times.
\item[S.6]
The first four stochastic elements (S.1)--(S.4) are mutually independent.
\end{enumerate}

The notation $x^i(N)$, $\sig^i(N)$ and $\eta^l(N)$
is henceforth abbreviated to $x^i$, $\sig^i$ and $\eta^l$.

The construction based on these primitives in presented momentarily,
but first it is helpful to clarify two points.
First,
the reason $\eta^l$ are called removal attempt times, not removal times,
is that it is possible that there are no particles in the system when one of these
times occurs, in which case no particle is actually removed
(see the 4-th bullet below).
Second, by the independence assumption,
a.s., no simultaneous introduction of particles (by injection or birth)
can occur after time $0$.
However, the introduction of a new particle and the removal of a particle
may occur simultaneously, such as when branching and removals
are coupled (see the last bullet below).

The term {\it leftmost particle}
refers to the particle whose position is lowest among the particles
in the configuration at a given time. Ties are broken according to some fixed
ordering of the labels. A particle is said to be {\it introduced} into the system
at a certain time if it is either injected (as in (S.2)) or born out of a branching event
(as is (S.4)).

The construction now proceeds in two steps.
First, once the initial space-time position $(x^i,\sig^i)$
of particle $i$ is determined, its {\it potential trajectory}, denoted
$\{X^i_t,t\in[\sig^i,\iy)\}$, is defined by
\begin{equation}\label{060}
X^i_t=\frX(x^i,\sig^i,B^i)(t),\qquad t\ge\sig^i.
\end{equation}
In the second step, the removal time $\tau^i$ of particle $i$ is determined
(where $\iy$ is possible),
and the {\it actual trajectory} the particle follows is obtained
by trimming the potential trajectory at $\tau^i$.

The particle configuration is defined on $(\eta^l,\eta^{l+1}]$
inductively for $l=0,1,2,\ldots$, where $\eta^0=0$.
The configuration at time $0$ is given by $X^i_0=x^i$ for $i=(j,0)$,
$j\in[N]$. This gives a well defined potential trajectory
of each of these particles on $[0,\iy)$.
Next, for $l\ge0$, given the configuration during $[0,\eta^l]$,
the construction during $(\eta^l,\eta^{l+1}]$ is described as follows.

During the time interval $(\eta^l,\eta^{l+1})$:
\begin{itemize}
\item
Each of the particles living at $\eta^l$ already has a well defined potential trajectory.
These particles live through the interval, with their actual trajectories given by
their potential trajectories.
\item
Each $i$ of the form $(j,0)$ with $\sig^i\in(\eta^l,\eta^{l+1})$
corresponds to an injection during this interval. This
determines the injection space-time location $(x^i,\sig^i)$ of a new particle,
and accordingly its potential trajectory for all $t\ge\sig^i$.
These particles live through the remainder interval.
\item
If a particle $i=(j,k)$ is alive when $\pi^i$ ticks, it gives birth
to a new particle at that space-time location.
The new particle gets the label $\hat i=(j,\hat k+1)$ where
$(j,\hat k)$ is the latest descendant of $(j,0)$ introduced prior to that time.
Again, this determines the potential trajectory,
and the particle lives through the remainder interval.
\end{itemize}

At the time $\eta^{l+1}$:

\begin{itemize}
\item
If there are no particles in the system (that is, the configuration
at time $\eta^l$ is empty), nothing happens. Otherwise:
\item
If no new particle is introduced at that time then the particle that is leftmost
at $\eta^{l+1}-$ is removed. If the index of this particle is $i$ then this
determines its removal time as $\tau^i=\eta^{l+1}$.
\item
If a particle is introduced at $\eta^{l+1}$ (by injection or branching),
the construction obeys the rule `introduce and then remove',
and this may cause the new particle to be removed immediately
(i.e., if the injection is to the left of all particles or the branching
particle is the leftmost).
\end{itemize}
For particles $i$ that never get removed, define $\tau^i=\iy$.
The lifetime of particle $i$ is given by $[\sig^i,\tau^i)$
(empty if $\sig^i=\tau^i$ or $\sig^i=\iy$) and
its actual trajectory is defined by $\{X^i_t:\,t\in[\sig^i,\tau^i)\}$.
This completes the construction of the particle system.

Some useful notation is as follows.
The set of particles initially in the system, injected and, respectively,
descendants of a root particle $i=(j,0)$, are denoted by
\[
\calS^{N,\rm init}=[N]\times\{0\},
\qquad
\calS^{N,\rm inj}=\{N+1,N+2,\ldots\}\times\{0\},
\qquad
\calS^{N,i}=\{(j,k):k\in\Z_+\}.
\]
The set of particles introduced by time $t$ is
\[
\calS^N_t=\{i\in\calS:\sig^i\le t\}.
\]
Those injected by time $t$ and, respectively, descendants of root particle $i$
introduced by time $t$, are denote by
\begin{align}\label{r20}
\calS^{N,\rm inj}_t=\calS^{N,\rm inj}\cap\calS^N_t,
\qquad
\calS^{N,i}_t=\calS^{N,i}\cap\calS^N_t.
\end{align}
Next, the configuration process is given by
\[
\xi^N_t(dx)
=\sum_{i\in\calS}\bdel_{X^i_t}(dx)1_{\{\sig^i\le t<\tau^i\}},
\]
and clearly its initial condition is
\[
\xi_0^N(dx)=\sum_{i\in\calS^{N,\rm init}}\bdel_{x^i}(dx).
\]
The injection and, respectively, removal
space-time locations are encoded by the random measures
\begin{equation}\label{r14}
\al^N(dx,dt)=\sum_{i\in\calS^{N,\rm inj}}\bdel_{(x^i,\sig^i)}(dx,dt),
\qquad
\beta^N(dx,dt)=\sum_{i\in\calS:\,\tau^i<\iy}\bdel_{(X^i_{\tau^i},\tau^i)}(dx,dt).
\end{equation}
Let $m^N_t=|\xi^N_t|$ denote the number of living particles at $t$.
Let $I^N_t=\#\calS^{N,\rm inj}_t$ be the number of injections by time $t$.
Then $I^N_t=\al^N(\R\times[0,t])$. Moreover,
let $J^N_t=\#\{l\ge1:\eta^l\le t\}$ be the number of removal attempts by $t$.
Note that the actual number of removals by $t$ is $\beta^N(\R\times[0,t])$.
Then $J^N_t=\beta^N(\R\times[0,t])$ holds on the event
$\{\inf_{s\le t}m^N_s\ge1\}$.

So far we have not made any assumption on the removal attempt times $\eta^l$.
We would like to cover the possibility that removals are coupled with
(some of the) injections or branching events, as well as
that they occur independently of each other. Hence we let
\begin{equation}\label{r8}
\calF^N_t=\sig\{\xi^N_0,\al^N(-\iy,x]\times[0,s],B^i_s,\pi^i_s,J^N_s:
x\in\R,s\in[0,t],i\in\calS\},
\end{equation}
and supplement (S.1)--(S.6) above with
\begin{itemize}
\item[S.7]
$B^i$ and $\hat\pi^i$ are $\{\calF^N_t\}$-martingales,
$i\in\calS$, where $\hat\pi^i(t)=\pi^i(t)-\kap t$.
\end{itemize}

\subsection{Macroscopic data}


\begin{definition}\label{def1}
An {\em admissible macroscopic data} is a
deterministic tuple $(\xi_0,\al,J)$
satisfying the following conditions, for some $\rho_0>0$.

i.
$\xi_0\in\calM_1(\R)$, $\al\in\calM_{+,\rm loc}(\R\times\R_+)$
and $J\in AC^\up(\R_+,\R_+)\cap C^{\rho_0}(\R_+,\R_+)$.

ii. Denote $I_t=\al(\R\times[0,t])$. Then one of the following holds:

\quad 1.
$\al(dx,dt)\sqsubset c\,dx\,dI_t$, some constant $c$. Moreover,
$I\in C^\up(\R_+,\R_+)\cap C^{\rho_0}(\R_+,\R_+)$.

\quad 2.
$I\in C^\up(\R_+,\R_+)\cap C^{\frac{1}{2}+\rho_0}(\R_+,\R_+)$.

iii. Let $m$ denote the solution to
\begin{equation}\label{055}
m_t=1+\kap\int_0^tm_sds+I_t-J_t,
\end{equation}
representing the total macroscopic mass. Then $\eps_0:=\inf_{t\in\R_+}m_t>0$.

\end{definition}

The stochastic elements of the model, namely
$\al^N$ and $\beta^N$ (respectively, $\xi^N_0$; $\xi^N$; $I^N$ and $J^N$)
are viewed as random variables taking values in $\calM_{+,\loc}(\R\times\R_+)$
(respectively, $\calM_1(\R)$; $D(\R_+,\calM_+(\R))$; $D(\R_+,\R)$).
The assumed structure of the stochastic primitives is as follows
(recall that the bar notation stands for normalization by $N$).

\begin{assumption}\label{assn1}
As $N\to\iy$, $(\bar\xi^N_0,\bar\al^N,\bar J^N)\to(\xi_0,\al,J)$ in the product topology,
in probability, where the latter is an admissible macroscopic data.
\end{assumption}

\begin{remark}\label{rem01}
i. ($N$-BBM as a special case).
In the $N$-BBM model, branching occurs at rate $1$ per particle,
and removal and branching
are coupled, so that the number of particles remains $N$ at all times.
Thus the birth/removal counting process, $J^N$, is Poisson of rate $N$. This model satisfies our assumptions:
Assumption \ref{assn1} is satisfied
with $\al=0$, $\kap=1$, $J_t=t$ and $(\frb,\frc)=(0,1)$, while
condition (S.7) can be seen to hold under the aforementioned
coupling.
Next, if one abandons the requirement
that births and removals are coupled, and instead assumes that
the removal count $J^N$ and the branching mechanism are mutually
independent, but $J^N$ is still Poisson of rate $N$, then all
assumptions are still valid with the same macroscopic data.
This latter version can be extended to allow variable removal rate
by letting $J\in AC^\uparrow(\R_+,\R_+)$ and assuming that $m_t$ defined by
\[
m_t=1+\int_0^tm_sds-J_t
\]
remains positive at all times. Letting $J^N$ be an inhomogeneous Poisson process
with instantaneous intensity $NJ'_t$ results in
a BBM with macroscopic mass $m_t$.

ii. (The model from \cite{de-masi-book}).
Similarly, the injection-selection model of \cite{de-masi-book}, mentioned
in the introduction, is closely related.
In this model the injections and removals are again coupled.
Consider $\kap=0$,
$\al^N$ a Poisson point process with intensity $N\pi(dx)c_0dt$, where $\pi$ is
any probability measure and $c_0>0$ a constant, and $J^N=I^N=\al^N(\R\times[0,t])$.
Then the case $\pi=\bdel_0$ gives the model from \cite{de-masi-book}
except the minor point
that the particles live in $\R$ rather than $\R_+$.
\end{remark}

\subsection{Weak FBP formulation and main result}\label{sec2.2}

First we recall the notion of second order parabolic equations with measure-valued r.h.s.\ \cite{ama01}.
Let $\fra=\frac{1}{2}\frc^2$ and
\begin{equation}\label{r9}
\calL\ph=\fra\pl_x^2\ph+\frb\pl_x\ph+\kap\ph,
\qquad
\calL^*u=\pl_x^2(\fra u)-\pl_x(\frb u)+\kap u.
\end{equation}
For $\xi_0\in\calM_1(\R)$, $\mu\in\calM_{\rm loc}(\R\times\R_+)$ consider the equation
\begin{equation}\label{017}
\begin{cases}
\pl_tu-\calL^*u=\mu,\\
u(\cdot,0)=\xi_0.
\end{cases}
\end{equation}
Let $q\in(1,\iy)$. A weak $L_q$-solution of \eqref{017} is a function
$u\in L_{1,\rm loc}(\R_+,L_q)$ satisfying
\begin{equation}\label{r4}
-\int_0^\iy(\pl_t\ph+\calL\ph,u)dt
=\int_\R\ph(\cdot,0)\,d\xi_0+\int_{\R\times\R_+}\ph\,d\mu
\end{equation}
for all $\ph\in C_c^\iy(\R\times\R_+)$.

Such problems were analyzed in \cite{ama01}.
In particular, \cite[Theorems 1 and 13 and Remarks 1 and 2(b)]{ama01} show that
for $1<q<\iy$, this problem possesses a unique weak $L_q$-solution,
independent of $q$
(note
that with the transformation $\bar\fra=\fra$, $\bar\frb=-\frb+\fra'$,
$\bar\kap=\kap-\frb'$, one has
the divergence form $\calL^* u=\pl_x(\bar\fra\pl_xu)+\bar\frb\pl_xu+\bar\kap u$,
as required in \cite[Remark 1(e)]{ama01};
also note that in \cite{ama01} the initial condition is accounted for
by substituting $\mu+\xi_0\otimes\bdel_0$ for $\mu$).
In what follows we thus use the term
weak solution to \eqref{017}, without reference to $q$.

We base on this notion the following problem formulation.
Let admissible data $(\xi_0,\al,J)$ be given. Consider the equation
\begin{equation}\label{018}
\begin{cases}
{\it (i)}\hspace{3.5em}
\pl_tu-\calL^*u=\alpha-\beta,\\
{\it (ii)}\hspace{3.1em}
\beta(U>0)=0 &\text{where} \quad \ds U(x,t)=\int_{-\iy}^x u(y,t)dy,\\
{\it (iii)}\hspace{2.7em}
\beta(\R\times[0,t])=J_t & \text{for}\quad t\in\R_+,\\
{\it (iv)}\hspace{2.8em}
u(\cdot,0)=\xi_0.
\end{cases}
\end{equation}
A solution $(u,\beta)$ to \eqref{018} is defined as a member of
$L_{1,\rm loc}(\R_+,L_q)\times\calM_{+,\rm loc}(\R\times\R_+)$
for some (equivalently, all) $q\in(1,\iy)$, such that $u$ is
an a.e.\ non-negative weak solution to \eqref{018}(i, iv),
and moreover, conditions \eqref{018}(ii, iii) hold.

For future reference, according to \eqref{r4}, a weak solution to \eqref{018}(i, iv)
is one for which
\begin{equation}\label{054}
-\int_0^\iy(\pl_t\ph+\calL\ph,u)dt
=\int_{\R}\ph(\cdot,0)d\xi_0+\int_{\R\times\R_+}\ph d\al
-\int_{\R\times\R_+}\ph d\beta,
\end{equation}
for $\ph\in C_c^\iy(\R\times\R_+)$.

\begin{theorem}\label{th2}
Let Assumptions \ref{assn0} and \ref{assn1} hold. Then
\\
i.
There exists a unique solution $(u,\beta)$ to \eqref{018}.
\\
ii.
There exists a version of $u$, again denoted $u$, such that setting
$\xi_t(dx)=u(x,t)dx$ gives $\xi\in C(\R_+,\calM_+(\R))$, and
$(\bar\xi^N,\bar\beta^N)\to(\xi,\beta)$ in probability as $N\to\iy$.
\end{theorem}

\begin{remark}\label{rem00}
i. The notion \eqref{018} can be seen as an extension of
a classical solution to a FBP. For example, consider
\eqref{402} and assume that $(u,\ell)$ is a classical solution,
with $\ell\in C(\R_+,\R)$ and $u$ appropriately smooth.
Recall that the injection-removal rate is set to $c_0$
in \eqref{018}.
Then a solution $(u,\beta)$ to \eqref{018} is obtained by
$\beta(dx,dt)=\bdel_{\ell_t}(dx)c_0dt$, as can be verified directly.

ii. More generally, one can construct a free boundary out of
any solution $(u,\beta)$ of \eqref{018}. If we take, as in Theorem \ref{th2},
$\xi_t(dx)=u(x,t)dx$, and disintegrate
$\beta$ as $\beta(dx,dt)=\beta_t(dx)dJ_t$, then two candidates are given by
\[
\ell_t=\inf\supp\, \xi_t,
\qquad
\tilde\ell_t=\sup\supp\,\beta_t,
\qquad
t\in(0,\iy).
\]
These may in general behave very irregularly.
For example, consider $\calL^*u=\pl_x^2u+u$ and $\al=0$,
and for $K\in\calB(\R_+)$, let $J_t=|[0,t]\cap K|$. (Note that $J_t$
satisfies our assumptions and thus $\xi,\beta,\ell$ and $\tilde\ell$ are well defined.)
Now, if $K=[0,1]\cup[2,\iy)$ then
during the time interval $(0,1)$, the free boundary, say $\ell$,
evolves continuously (by results of \cite{ber19}), at time 1
it jumps to $-\iy$ (as there is no absorption of mass),
and at time 2 it comes back from $-\iy$. It is clear that
this behavior can be made far more complicated by picking
other choices of $K$, and for a general Borel set $K$
one does not expect any regularity of $\ell$ beyond Borel measurability.

Under such circumstances,
the notion of a free boundary does not seem
to be particularly useful from an analytic point of view
(e.g.\ since boundary conditions in initial boundary value problems
are typically specified on the boundary of an open set).
Also, a Feynman-Kac representation of
the component $u$ of the solution in terms of $\ell$,
which involves the hitting time of BM to $\ell$, requires
some regularity of $\ell$ (for example, piecewise continuity).
Feynman-Kac representation of $u$ in terms of $\ell$
has served as a key tool in earlier work on the subject
(see details in \S \ref{sec30++}), but, for the reasons mentioned, is not available here.
\end{remark}


\section{Uniqueness via barriers}\label{sec3}
\beginsec

In this section we prove the following result.
\begin{theorem}\label{th1}
Let Assumption \ref{assn0} hold and let admissible data $(\xi_0,\al,J)$
be given.
Then there exists at most one solution $(u,\beta)$ to \eqref{018}.
\end{theorem}

The proof is based on the construction of barriers, which are shown to constitute upper
and lower bounds to any solution, in the sense of mass transport inequalities. This section is structured as follows.
Essential tools are developed in \S\ref{sec30} and \S\ref{sec30+},
where the former provides so called mild solutions,
and the latter introduces operators required for the construction,
and studies some of their properties. \S\ref{sec30++} gives a brief description of
the construction of barriers in earlier work and of the difference to that
in the current work, as well as a sketch of the proof of uniqueness based on them.
The upper and lower barriers are constructed in \S\ref{sec31} and \S\ref{sec32}
respectively. In \S\ref{sec33} it is shown that the barriers
can be made close to each other, and the proof is completed.

\subsection{Preliminary lemmas}\label{sec30}

The backward Kolmogorov equation associated with the diffusion \eqref{407}
is given by $\pl_tu=\calL_1u$ with $\calL_1u=\fra\pl_x^2u+\frb\pl_xu$.
Denote by $\frp_t(x,y)$ the fundamental solution of this equation.

\begin{lemma}\label{lem8}
Given $T$ there exist constants $\tilde c_1,\tilde c_2, \hat c_1, \hat c_2>0$ such that
for $t\in(0,T]$ and $x,y\in\R$,
\begin{equation}\label{91}
\tilde c_1t^{-1/2}e^{-\tilde c_2(x-y)^2t^{-1}}
\le \frp_t(x,y)\le \hat c_1t^{-1/2}e^{-\hat c_2(x-y)^2t^{-1}}.
\end{equation}
\end{lemma}
Whereas the upper bound will be used many times,
the lower bound is needed only to make the following
statement (used in Lemma \ref{lem2}):
There exists a constant $c_*>0$ such that
\begin{equation}\label{92}
\int_{-\iy}^x\frp_t(x,y)dy\ge c_*,\qquad t\in(0,1], \, x\in\R.
\end{equation}
\begin{proof}
For $T=1$ these bounds follow from \cite[Theorems 4.4.6 and 4.4.12]{stroock}.
To verify the assumptions, note that one can write $\calL_1$
in the form $\calL_1u=\pl_x(\fra\pl_xu)+\tilde\frb\fra\pl_xu$
by setting $\tilde\frb=(\frb-\fra')\fra^{-1}$. The boundedness of $\tilde\frb$ follows
from the assumed boundedness of $\fra^{-1}$, $\fra'$ and $\frb$.
For $T>1$, apply the scaling property of $\frp$ as in \cite[Remark 4.1.5]{stroock}
(with constants depending on $T$).
\end{proof}

Denote
\begin{equation}\label{93}
\frs_t(x,y)=e^{\kap t}\frp_t(x,y).
\end{equation}
For $u\in L_1(\R)$ denote
\begin{equation}\label{r1}
S_tu(y)=\int_\R \frs_t(x,y)u(x)dx,
\end{equation}
and with a slight abuse of notation, use the same symbol for $\xi\in\calM_+(\R)$, namely
\begin{equation}\label{r2}
S_t\xi(y)=\int_\R\frs_t(x,y)\xi(dx).
\end{equation}
For $\gamma\in\calM_{+,\loc}(\R\times\R_+)$ and $0\le \tau<t$, denote
\[
S*\gamma(y,t)=\int_{\R\times[0,t]}\frs_{t-s}(x,y)\gamma(dx,ds),
\]
\[
S*\gamma(y,t;\tau)=\int_{\R\times[\tau,t]}\frs_{t-s}(x,y)\gamma(dx,ds).
\]

\begin{lemma}\label{lem1}
i. Let $\gam\in\calM_{+,\rm loc}(\R\times\R_+)$ be such that
$\gam(\R\times[0,\cdot])\in C^{\rho_0}(\R_+,\R_+)$ for some
$\rho_0>0$. Let $v(y,t)=S_t\xi_0(y)+S*\gam(y,t)$.
Then, for $q\in(1,\iy)$, $v\in L_{1,\rm loc}(\R_+,L_q)$.

ii. Let $(\tilde u,\beta)$ be a solution to \eqref{018}.
If $u$ is a version of $\tilde u$ then $(u,\beta)$ is also a solution.
Moreover, $\tilde u$ has a version given by
\begin{equation}\label{012-}
u(y,t)=S_t\xi_0(y)+S*\al(y,t)-S*\beta(y,t).
\end{equation}

iii. One has $\|u(\cdot,t)\|_1=m_t$, $t>0$. Moreover,
$v,u\in L_{\iy,\loc}((0,\iy),L_\iy)$.

iv. One has, for $0<\tau<t$,
\begin{equation}\label{012}
u(y,t)=S_{t-\tau}u(\cdot,\tau)(y)+S*\al(y,t;\tau)-S*\beta(y,t;\tau).
\end{equation}
\end{lemma}

\begin{remark}\label{rem2}
In what follows,
by a solution to \eqref{018} we will mean the version given
by \eqref{012-} unless stated otherwise. In view of Lemma \ref{lem1},
there is no loss of generality in doing so when proving the uniqueness result.
\end{remark}

\begin{proof}
i.
Fix $T$. In this proof,
$c$ denotes a constant not depending on $x,y$ and $t\in(0,T]$
whose value may change from one expression to another.
By Lemma \ref{lem8} and \eqref{93}, it is easy to see that
$\|\frs_t(x,\cdot)\|_2\le ct^{-1/2}$, $t\in(0,T]$. By Minkowski's integral inequality
it follows that $\|S_t\xi_0\|_2\le ct^{-1/2}$.

Next, let $q\in(1,\iy)$. Then by
Lemma \ref{lem8} and \eqref{93}, for $t\in(0,T]$,
$\|\frs_t(x,\cdot)\|_q\le c t^{-Q}$ where $Q=(q-1)/(2q)$. By Minkowski's integral inequality,
\[
\|S*\gamma(\cdot,t)\|_q\le c\int_{\R\times[0,t]}(t-s)^{-Q}\gamma(dx,ds)
=c\int_{[0,t]}(t-s)^{-Q}dK_s,
\]
where $K_t=\gam(\R\times[0,t])$.
By monotone convergence, the last integral is the limit as $\eps\downarrow0$ of
\begin{align}
\notag
\int_{[0,t-\eps]}(t-s)^{-Q}dK_s
&=
K_{t-\eps}\eps^{-Q}-Q\int_0^{t-\eps}K_s(t-s)^{-1-Q}ds
\\
&=(K_{t-\eps}-K_t)\eps^{-Q}
+Q\int_0^{t-\eps}(K_t-K_s)(t-s)^{-1-Q}ds+K_tt^{-Q}
\notag
\\
&\le c\int_0^t(t-s)^{\rho_0-1-Q}ds+K_tt^{-Q}.
\label{r5}
\end{align}
If we choose $q>1$ sufficiently small then $\rho_0-1-Q\in(-1,0)$,
and the above integral is bounded by $c$ for $t\le T$. Moreover,
$Q\in(0,1)$, and we obtain that $\|S*\gam(\cdot,t)\|_q$ is
integrable over $[0,T]$.

The two terms defining $v$ satisfy \eqref{r4} with data
$(\xi_0,0)$ and $(0,\gamma)$, respectively, by elementary integration
by parts (see e.g. \cite[Theorem 4.6]{folland95} for a similar calculation).
The estimates above on these two terms show that they are, respectively,
members of $L_{1,\rm loc}(\R_+,L_2)$,
and $L_{1,\rm loc}(\R_+,L_q)$ for $q$ close to $1$.
Hence each is a weak $L_q$-solution of the corresponding equation
for some $q\in(1,\iy)$.
By the results of \cite{ama01} discussed following \eqref{r4},
they must therefore be the unique
weak $L_q$-solution, for all $q\in(1,\iy)$. This proves the assertion.

ii.
For the first assertion we must show that the complementary condition
\eqref{018}(ii) is preserved by changing $\tilde u$ to $u$.
If $U$ is as in \eqref{018}(ii) and $\tilde U$ is defined similarly,
then there is a set $A\subset\R_+$ of full Lebesgue measure
such that $U(x,t)=\tilde U(x,t)$ for $(x,t)\in\R\times A$.
Owing to the assumption that $J$ is absolutely continuous,
$\beta$ does not charge $\R\times A^c$. This shows that
$\beta(U>0)$ holds if and only if $\beta(\tilde U>0)$, proving the assertion.

For the second assertion,
the arguments given above in (i) show that the three terms
on the right of \eqref{012-} are, for every $q\in(1,\iy)$,
weak $L_q$-solutions of \eqref{017}
for $\mu=\xi_0\otimes\bdel_0$, $\al$ and, respectively, $-\beta$.
By linearity in $\mu$ of weak solutions of \eqref{017} \cite[Theorem 1]{ama01},
it follows that $u$ defined in \eqref{012-} is a weak $L_q$-solution
of \eqref{018}(i) corresponding to data $(\xi_0,\al,\beta)$.
In particular, $u$ must be a version of $\tilde u$ by uniqueness
of solutions to \eqref{017}.

iii.
To calculate $\tilde m_t:=\|u(\cdot,t)\|_1$, note that
$\|\frs_t(x,\cdot)\|_1=e^{\kap t}$.
Hence by \eqref{012-},
\[
\tilde m_t=e^{\kap t}+\int_0^te^{\kap(t-s)}dI_s-\int_0^te^{\kap(t-s)}dJ_s,
\]
which solves \eqref{055}, and by uniqueness, equals $m_t$.

As for the estimate on $\|u(\cdot,t)\|_\iy$, by positivity
we only need to estimate the first two terms in \eqref{012-}.
Directly from Lemma \ref{lem8}, the sum of these two terms is bounded as follows:
\[
\frs_t(x,y)\le ct^{-1/2}+cI_t
\qquad\text{and}\qquad
ct^{-1/2}+c\int_0^t(t-s)^{-1/2}dI_s,
\]
under Assumption \ref{assn1}(ii.1) and, respectively, (ii.2).
The former expression is locally bounded for $t$ away from $0$,
as required. As for the latter, a calculation as in \eqref{r5},
replacing $(Q,\rho_0,K_t)$ by $(\frac{1}{2},\frac{1}{2}+\rho_0,I_t)$,
shows that this expression is also locally bounded for $t$ away from $0$.

iv.
Finally, \eqref{012} follows from \eqref{012-} upon using the Chapman-Kolmogorov equation
$\int \frp_{t-\tau}(x,y)\frp_\tau(y,z)dy=\frp_t(x,z)$.
\end{proof}

\subsection{Mass transport inequalities}\label{sec30+}

In this section, several elementary facts about mass transport inequalities are borrowed
from \cite{de-masi-book}, and some are developed further.
On $L_1(\R,\R_+)$, define the relation $u\LE v$ as
\[
u[r,\iy)\le v[r,\iy) \quad \text{for all } r\in\R,
\]
and the relation $u\LE v\MOD \ell$, for $\ell\ge0$, as
\[
u[r,\iy)\le v[r,\iy)+\ell \quad \text{for all } r\in\R.
\]
For $\xi,\zeta\in\calM_+(\R)$, define $\xi\LE\zeta$ and $\xi\LE\zeta\MOD\ell$ analogously.

For $\del>0$, the `cut' operator $C_\del$ acts on $H_\del=\{u\in L_1(\R,\R_+):
\|u\|_1>\del\}$
by cutting mass of size $\del>0$ on the left. That is, for $u\in H_\del$,
\[
\La_\del(u)=\inf\{x\in\R:u(-\iy,x] \ge \del\}
\]
and
\[
C_\del u(x)=u(x)1_{[\La_\del(u),\iy)}(x).
\]
When $\del=0$ set $C_\del=\id$, the identity map.
Also, denote $\vhat C_\del={\tt id}-C_\del$. We also use an operator
that cuts out a mass of size $\del$ lying between $\La_\Del$ and $\La_{\Del+\del}$.
More precisely, given $\Del>0$ and $\del\ge0$, $\DEL$ will always denote
the pair $(\Del,\del)$.
Then the operator $C_\DEL$ acts on $H_{\Del+\del}$
as
\[
C_{\DEL}u(x)=C_{\Del,\del}u(x)=u(x)1_{\{(-\iy,\La_{\Del}(u)]\cup(\La_{\Del+\del}(u),\iy)\}}(x).
\]
Set $\vhat C_\DEL=\id-C_\DEL$.

\begin{lemma}\label{lem3}
Let $\del\ge0$ and assume $u,v\in H_\del$. Let $\ell\ge0$.

i. If $u\LE v\MOD \ell$ and $\|u\|_1=\|v\|_1$ then $S_\del u\LE S_\del v\MOD e^{\kap\del}\ell$.

ii. If $u\LE v\MOD \ell$ then $C_\del u\LE C_\del v\MOD \ell$.

iii.
If $w\in L_1(\R,\R_+)$ is such that $\|w\|_1=\del$
and $u-w\ge0$ then $u-w\LE C_\del u$.

Next let $\Del>0$ and assume $u,v\in H_{\Del+\del}$.

iv. If $u\LE v \MOD \ell$ then $C_{\Del,\del}u\LE C_{\Del,\del}v \MOD \ell$.

v. If $0<\vhat\Del\le\Del$ then $C_{\Del,\del}u\LE C_{\vhat\Del,\del}u$.

vi. If $\Del'\ge\Del$, $u\LE v\MOD \Del'$ and $\|u\|_1=\|v\|_1$
then $C_\del u\LE C_{\Del,\del}v\MOD\Del'$.
\end{lemma}

\begin{proof}
i.
Step 1. Consider the case where $\kap=0$ and
$\|u\|_1=\|v\|_1=1$.
In this case $u,v$ are densities of probability measures.
Without loss of generality (w.l.o.g.) assume $\ell<1$. Fix $r_0$ be such that
$\int_{r_0}^\iy u=\ell$ (where $r_0=\iy$ if $\ell=0$).
Let $\tilde u=u1_{\{\cdot\le r_0\}}$ be a density that integrates to $1-\ell$.
Consider the probability measure on $[-\iy,\iy)$, denoted $\tilde U$,
having mass $\ell$ at $-\iy$ and density $\tilde u$ on $\R$.
Let $V$ be the probability measure with density $v$.
Then one has $\tilde U[r,\iy)\le V[r,\iy)$ for all $r$.
Therefore the exists a coupling $(\tilde X_0,Y_0)$
having marginal distributions $\tilde U$ and $V$, respectively,
such that $\tilde X_0\le Y_0$ a.s. Denote by $E$ the event
$\{\tilde X_0>-\iy\}$.

Consider a coupling of two processes $\tilde X$ and $Y$,
constructed using a BM $B$ independent of $(\tilde X_0,Y_0)$.
Namely, on the event $E$ let $\tilde X_t$ be the unique strong solution to
\[
\tilde X_t=\tilde X_0+\int_0^t\frb(\tilde X_s)ds+\int_0^t\frc(\tilde X_s)dB_s.
\]
$\tilde X_t$ need not be defined on $E^c$.
Similarly, define $Y_t$ (on all of $\Om$)
as the solution to this SDE with initial condition $Y_0$.
Then $\tilde X_t\le Y_t$ for all $t$ holds a.s.\ on $E$. This gives
\[
\PP(\tilde X_\del>r)\le \PP(E\cap\{Y_\del>r\})
\le \PP(Y_\del>r)=S_\del v(r,\iy).
\]
Next, let $X_0=\tilde X_0$ on $E$,
and let its conditional law given $E^c$ be given by
the density $\ell^{-1}u1_{\{\cdot>r_0\}}$ (which need not be defined
in the case $\ell=0$). Then $X_0$ has $u$ as its density.
Again, assume w.l.o.g.\ that $B$ is independent of $X_0$,
and let $X_t$ be defined (on all of $\Om$)
by following the same SDE with $X_0$
as an initial condition. Then the density of $X_\del$ is given
by $S_\del u$, and $X_t=\tilde X_t$ on $E$. Thus for any $r\in\R$,
\[
S_\del u(r,\iy)=\PP(X_\del>r)
\le\PP(E\cap\{\tilde X_\del>r\})+\PP(E^c)\le S_\del v(r,\iy)+\ell.
\]

Step 2. If $\|u\|_1=\|v\|_1=c$ then $u/c$ and $v/c$ are probability
densities and $u/c\LE v/c$ mod $\ell/c$. This gives, by Step 1,
$S_\del u/c\LE S_\del v/c$ mod $\ell/c$. The claim follows
on multiplying by $c$.

Step 3. Finally, when $\kap>0$, the claim follows from Step 2 after
multiplying by $e^{\kap\del}$ and using \eqref{93}.

ii.
Let $a=\La_\del(u)$ and $b=\La_\del(v)$.
For $r\ge a\vee b$, clearly
$(C_\del u)[r,\iy)=u[r,\iy)\le v[r,\iy)+\ell=(C_\del v)[r,\iy)+\ell$.
For $r< a\vee b$, consider two cases.

Case 1: $a\le b$ and $r< b$. Then
\[
(C_\del u)[r,\iy)\le\|C_\del u\|_1=\|u\|_1-\del
\le\|v\|_1-\del+\ell=\|C_\del v\|_1+\ell=C_\del[r,\iy)+\ell.
\]

Case 2: $a>b$ and $r< a$. Then
\[
(C_\del u)[r,\iy)=u[a,\iy)\le u[r\vee b,\iy)\le v[r\vee b,\iy)+\ell=(C_\del v)[r,\iy)+\ell.
\]

iii.
We have $C_\del u=u1_{[c,\iy)}$, where $u(-\iy,c]=\del$.
Consider $r\le c$. Because $u-w$ is nonnegative,
\[
(u-w)[r,\iy)\le \|u-w\|_1=\|u\|_1-\|w\|_1=\|u\|_1-\del=
(C_\del u)[r,\iy).
\]
Next, consider $r>c$. Then $(u-w)[r,\iy)\le u[r,\iy)$
whereas $(C_\del u)[r,\iy)=u[r,\iy)$.

iv. We have $u[r,\iy)\le v[r,\iy)+\ell$ for all $r$. We must show that
$\vhat u[r,\iy)\le\vhat v[r,\iy)+\ell$ for all $r$, where
\[
\vhat u=u1_{(-\iy,a)\cup(b,\iy)},
\text{ $u(-\iy,a)=\Del$, $u(a,b)=\del$},
\]
\[
\vhat v=v1_{(-\iy,\bar a)\cup(\bar b,\iy)},
\text{ $v(-\iy,\bar a)=\Del$, $v(\bar a,\bar b)=\del$}.
\]
We split into four cases.

Case 1. $r\le a$:
\[
\vhat u[r,\iy)=u[r,\iy)-\del\le v[r,\iy)-\del + \ell \le\vhat v[r,\iy) + \ell.
\]

Case 2. $r\ge b\vee\bar b$:
\[
\vhat u[r,\iy)=u[r,\iy)\le v[r,\iy) + \ell =\vhat v[r,\iy) + \ell.
\]

Case 3. $\bar b\le b$ and $a\le r<b$:
\[
\vhat u[r,\iy)=u[b,\iy)\le v[b,\iy)+\ell =\vhat v[b,\iy)+\ell \le \vhat v[r,\iy)+\ell.
\]

Case 4. $b<\bar b$ and $a\le r<\bar b$:
Note that
\[
\vhat u[b,\iy)=\|u\|_1-(\Del+\del),
\qquad
\vhat v[\bar b,\iy)=\|v\|_1-(\Del+\del).
\]
Moreover, $\|u\|_1\le\|v\|_1+\ell$. Hence $\vhat u[b,\iy)\le\vhat v[\bar b,\iy)+\ell$.
Therefore
\[
\vhat u[r,\iy)\le\vhat u[b,\iy)\le\vhat v[\bar b,\iy)+\ell \le\vhat v[r,\iy)+\ell.
\]

v.
Note that
\[
C_{\vhat\Del,\del}u= C_\del v+z,\qquad
C_{\Del,\del}u=C_{\Del-\vhat\Del,\del}v+z,
\]
where
\[
v=C_{\vhat\Del}u,\qquad z=\vhat C_{\vhat\Del}u.
\]
Hence it suffices to prove $C_{\Del-\vhat\Del,\del}v\LE C_\del v$.
To this end, note that $C_{\Del-\vhat\Del,\del}v=v-w\in L_1(\R,\R_+)$,
and $\|w\|_1=\del$. Therefore we can use part (iii) of the lemma,
by which $v-w\LE C_\del v$.
This completes the proof.

vi.
First, let us compare $u$ to $C_\Del v$.
If $r<\La_\Del(v)$ then
\[
u[r,\iy)\le\|u\|_1=\|C_\Del v\|_1+\Del
=C_\Del v[r,\iy)+\Del\le C_\Del v[r,\iy)+\Del'.
\]
If $r\ge\La_\Del(v)$ then
\[
u[r,\iy)\le v[r,\iy)+\Del'=C_\Del v[r,\iy)+\Del'.
\]
This shows that $u\LE C_\Del v \MOD \Del'$.
By part (ii) of the lemma,
$C_\del u\LE C_{\Del+\del}v \MOD \Del'$.
Finally, $C_{\Del+\del}v\le C_{\Del,\del}v$ pointwise, hence
$C_{\Del+\del}v\LE C_{\Del,\del}v$. This proves the claim.
\end{proof}

\subsection{On the barrier method}
\label{sec30++}

We refer to \cite{de-masi-book} for an exposition of the use of barriers
for proving uniqueness of solutions to FBP.
Here we briefly describe the main idea behind their use in earlier work,
and how their structure differs from this line in the present paper.

\vspace{.5em}

\noi{\it Barriers in earlier work.}
For simplicity, consider the $N$-BBM setting of \cite{de-masi-nbbm}.
Here, $\al=0$, $\kap=1$, $J_t=t$.
It is proved in \cite[Theorem 5]{de-masi-nbbm}
that for any classical solution $(u,\ell)$ to the FBP, $u$ is sandwiched between
lower and upper barriers in the sense that
\begin{equation}\label{105}
(S_\del C_{1-e^{-\del}})^n u_0\LE u(\cdot,n\del)
\LE (C_{e^\del-1}S_\del)^nu_0
\end{equation}
for all $\del>0$ and $n$. Uniqueness is then deduced by
showing that there can be at most one element that, for all $\del$ and $n$,
separates the lower from upper barriers.
(A similar argument appeared, for example,
in \cite[Theorem 3.14]{de-masi-book}
for relaxed FBP solutions and in \cite[Lemma 6.2]{ber19} for solutions
to the corresponding obstacle problem).
The proof of \eqref{105} is based on a Feynman-Kac representation of
$u$ in terms of $\ell$ \cite[eq.\ (6)]{de-masi-nbbm},
where the latter is assumed to be continuous.
(See, respectively, \cite[Proposition 8.1]{de-masi-book}
and \cite[Proposition 3.1]{ber19}).

Next, if $J$ is more general it is natural to expect, by analogy to \eqref{105}, that
\begin{equation}\label{106}
\prod_{i=0}^{n-1}(S_\del C_{\hat j_i(\del)}) u_0\LE u(\cdot,n\del)
\LE \prod_{i=0}^{n-1}(C_{j_i(\del)}S_\del)u_0
\end{equation}
where $\hat j_i(\del)$ and $j_i(\del)$ depend on $J$,
and reflect the fact that the mass removal rate varies.
However, as mentioned in Remark \ref{rem00}, a Feynman-Kac representation
for the density $u$ of \eqref{018} is missing in the generality of Theorem \ref{th2},
and therefore the argument given in the above references
requires adaptation. We can recover the
second inequality in \eqref{106} without this tool, but the first seems harder.
We construct an alternative
lower barrier for which we can prove the lower bound, by working with the operator
$C_{\Del,\delta}$ rather than $C_\del$.
Here our treatment considerably deviates from the above line.

\vspace{.5em}

\begin{figure}
\includegraphics[width=30em]{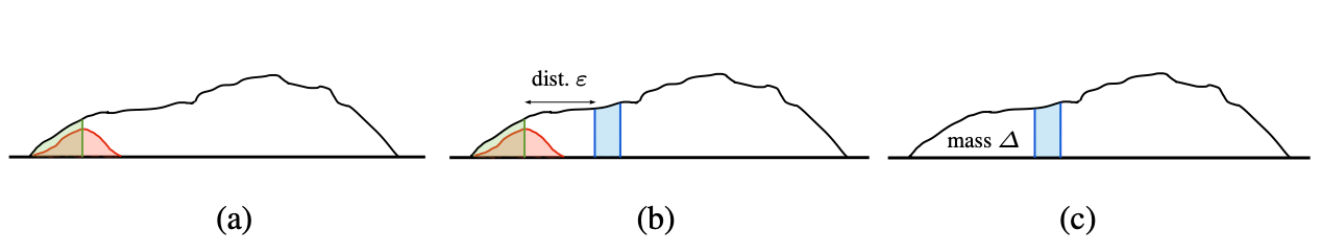}
\caption{(a) The solution is given by $v$ (black) minus $h$ (red);
upper barrier is obtained by removing the leftmost mass of size $j_0=\|h\|_1$
from $v$ (green).
(b) Removing mass of size $j_0$ (blue), located $\eps$ away from the mass in green.
(c) Removing instead mass of size $j_0$ (blue) after leaving mass of size $\Del$ on the left,
to obtain a lower barrier which can be iterated.
}
\label{fig-a}
\end{figure}

\noi{\it Sketch of construction of barriers in this work.}
The construction of the barriers and the proof that they form bounds
analogous to \eqref{105} and \eqref{106} appear
in the next two subsections. Here we sketch the main idea.
Let $(u,\beta)$ be a solution to \eqref{018}.
By equation \eqref{012-}, for $\del>0$, $u(\cdot,\del)=v-h$,
where
\[
v=S_\del\xi_0+S*\al(\cdot,\del),\qquad h=S*\beta(\cdot,\del).
\]
Let $j_0=\|h\|_1$. Then the nonnegative function $u(\cdot,\del)$
is obtained by removing
from $v$ the mass of size $j_0$ distributed according to $h$.
If instead one removes from $v$ the leftmost mass of size $j_0$,
as shown in Figure \ref{fig-a}(a),
then the resulting function $C_{j_0}v$ satisfies $u=v-h\LE C_{j_0}v$.
Next, by Lemma \ref{lem3}, both $S_\del$ and $C_{j_0}$ preserve the order,
and the argument can be iterated, providing an upper barrier at times $n\del$ for all $n$. In the case $\al=0$, this barrier takes
the form of the r.h.s.\ of \eqref{106}.

Next we sketch the lower barriers.
Figure \ref{fig-a}(b) shows in blue a mass of size $j_0$ located $\eps$ away from
the leftmost mass of size $\del$. If $\eps$ is fixed while $\del$ and $j_0$ are sufficiently
small, then one can show that most of the mass of $h$ (red) is to the left of
the mass in blue. Removing from $v$ the mass marked in blue
thus gives a lower barrier for $u$, up to an error term which can be made small.
While this is a valid statement, it is not useful for us because the operator that cuts
away the mass marked in blue does not preserve the order,
and as a result the inequality cannot be iterated. However, one can
use instead the operator $C_{\Del,j_0}$ that leaves
mass $\Del$ on the left and then cuts away mass $j_0$,
as shown in blue in Figure \ref{fig-a}(c). With an a priori bound on the density,
the statement regarding negligible mass in red reaching the mass in blue
is still valid here. Since, by Lemma \ref{lem3}, this operator preserves the order,
the argument can be iterated. As we will show,
the resulting error term can be controlled.

\subsection{Upper barriers}\label{sec31}

We now give the precise definitions.
The upper barriers are defined using a removal mechanism
that operates at times $n\del$, $n\in\N$, for $\del>0$ fixed.
Consider the time interval $[(n-1)\del,n\del]$.
The mass injected during this interval adds the term
$S *\al(\cdot,n\del;(n-1)\del)$ to the density at $n\del$. Hence let
the `paste' operator be defined, for $u\in L_1(\R,\R_+)$, by
\[
P^{(\del)}_nu=u+S*\al(\cdot,n\del;(n-1)\del).
\]
The mass removed during the said interval
is of size $J_{n\del}-J_{(n-1)\del}$.
However, more relevant is the size this mass would grow to be
had it not been removed, namely
\[
j_n(\del):=\int_{[(n-1)\del,n\del]}e^{\kap(n\del-s)}dJ_s.
\]
Accordingly, let
$C^{(\del)}_n=C_{j_n(\del)}$.

The upper barriers are defined for each $\del>0$ and $n\in\N$
by setting $u^{(\del,+)}_0=\xi_0$ and
\[
u^{(\del,+)}_{n\del}= C^{(\del)}_n P^{(\del)}_n S_\del u^{(\del,+)}_{(n-1)\del},\qquad n\in\N.
\]
Note that for $n=1$ and $n\ge2$, the function
$S_\del u^{(\del,+)}_{(n-1)\del}$ above is defined via \eqref{r2}
and, respectively, \eqref{r1}.
For the barriers to be well defined one must have
\begin{equation}\label{061}
P^{(\del)}_nS_\del u^{(\del,+)}_{(n-1)\del}\in H_{j_n(\del)}.
\end{equation}
We sometimes use the notation $u_t$ for $u(\cdot,t)$
as we do in the following.

\begin{proposition}\label{prop1}
Fix $\del>0$. Then \eqref{061} holds for all $n\in\N$, and consequently
the upper barriers are well defined. Moreover,
let $(u,\beta)$ be a solution to \eqref{018}.
Then for $n\in\N$,
\[
u_{n\del}\LE u^{(\del,+)}_{n\del}.
\]
\end{proposition}

\begin{proof}
Recall that by Assumption \ref{assn1}, $m_t>0$ for all $t$,
where $m_t$ is given by \eqref{055}. The $L_1$ norm of $u^{(\del,+)}$ satisfies
\[
\|u^{(\del,+)}_{n\del}\|_1=\|u^{(\del,+)}_{(n-1)\del}\|_1e^{\kap\del}
+\int_{[(n-1)\del,n\del]}e^{\kap(n\del-s)}dI_s
-\int_{[(n-1)\del,n\del]}e^{\kap(n\del-s)}dJ_s,
\]
so long as the r.h.s.\ above is positive.
By induction on $n$, this expression gives $\|u^{(\del,+)}_{n\del}\|_1=m_{n\del}>0$,
completing the proof of the first assertion.
Note by a simple induction argument, the definition of $j_n(\del)$
and Lemma \ref{lem1}, that $\|u_{n\del}\|_1=\|u_{n\del}^{(\del,+)}\|_1$.

The main claim is also proved by induction.
For $n-1\ge1$, assume $f\LE g$ where $f=u_{(n-1)\del}$ and $g=u^{(\del,+)}_{(n-1)\del}$; for $n-1=0$, $f=g=\xi_0$.
Write $C$, $P$ and $S$ for $C^{(\del)}_n$, $P^{(\del)}_n$ and $S_\del$, resp.
We have $u^{(\del,+)}_{n\del}=CPSg$.
Moreover, by Lemma \ref{lem1},
$u_{n\del}=PSf-h$, where
\[
h(y)=\int_{\R\times[(n-1)\del,n\del]}\frs_{n\del-s}(x,y)\beta(dx,ds).
\]
Hence the proof will be complete once $PSf-h\LE CPSg$ is shown.

By Lemma \ref{lem3}(i,ii), $C$ preserves the order $\LE$
and one has $S\tilde u\LE S\tilde v$ whenever $\tilde u\LE\tilde v$ and
$\|\tilde u\|_1=\|\tilde v\|_1$.
It is trivial that this is also true for $P$.
Denote $w=PSf$. Suppose one shows $w-h\LE Cw$. Then
\[
w-h\LE Cw=CPSf\LE CPSg,
\]
where in the last inequality one uses $Sf=Sg$ for $n-1=0$
and $\|f\|_1=\|g\|_1$ for $n-1\ge1$.
Thus the proof would be complete.

It thus suffices to show $w-h\LE Cw$.
The function $w-h$ is nonnegative
(as required by the definition of a solution)
and $\int h=j_n(\del)$. Hence $w-h\LE Cw$ by Lemma \ref{lem3}(iii),
and the proof is complete.
\end{proof}

\subsection{Lower barriers}\label{sec32}

The lower barriers are defined for $\DEL=(\Del,\del)\in(0,\iy)^2$
and $n\in\N$ as follows.
Let
\[
C^{(\DEL)}_n=C_{\Del,j_n(\del)}.
\]
Set $u^{(\DEL,-)}_0=\xi_0$ and $\ell_{0,\DEL}=0$, and for $n\in\N$,
\begin{align}\label{016}
u^{(\DEL,-)}_{n\del}
=C^{(\DEL)}_nP^{(\del)}_nS_\del u^{(\DEL,-)}_{(n-1)\del},
\qquad
\ell_{n,\DEL}=e^{\kap\del}\ell_{n-1,\DEL}+\begin{cases}
j_n(\del) & \text{if } (n-1)\del<t_0,\\
e^{-\Del^5/\del}j_n(\del)& \text{if } (n-1)\del\ge t_0,
\end{cases}
\end{align}
where $t_0>0$ is fixed.
Once again, for the definition to be valid, one must assure that for all $n\in\N$,
\begin{equation}\label{062}
P^{(\del)}_nS_\del u^{(\DEL,-)}_{(n-1)\del}
\in H_{\Del+j_n(\del)}.
\end{equation}

\begin{proposition}\label{prop2}
For $\Del\in(0,\eps_0)$ and $\del>0$, \eqref{062} holds for all $n$
and consequently the lower barriers are well defined.
Moreover, let $0<t_0<T$ be given.
Then there exists $\Del_0\in(0,\eps_0)$ such that for every $\Del\in(0,\Del_0)$
there exists $\del_0>0$ such that for $\del\in(0,\del_0)$
and $n\in\N$ satisfying $n\del\le T$ one has
\[
u^{(\DEL,-)}_{n\del}\LE u_{n\del} \MOD \ell_{n,\DEL}
\]
whenever $(u,\beta)$ is a solution to \eqref{018}.
Furthermore, for $n\in\N$, $n\del\le T$ one has
\begin{equation}\label{r3}
\ell_{n,\DEL}\le e^{\kap(T+\del)}(J_{t_0+\del}+e^{-\Del^5/\del}J_T).
\end{equation}
\end{proposition}

Given $\gamma\in\calM_{+,\loc}(\R\times\R_+)$ and $[t_1,t_2]\subset\R_+$,
the supremum of the support of the measure
$\tilde\gamma(\cdot)=\gamma(\cdot\times[t_1,t_2])$
is denoted by $\rho^*(\gamma;[t_1,t_2])$.
Recall $c_*$ from \eqref{92} and denote $c^*=2/c_*$.

\begin{lemma}\label{lem2}
Given a solution $(u,\beta)$ to \eqref{018}, $\del>0$ and $n\in\N$, let
\[
\rho_{n,\del}=\rho^*(\beta;[(n-1)\del,n\del]).
\]
Then for $n\ge2$,
\[
\rho_{n,\del}\le b(n,\del,u_{(n-1)\del}):=\La_{c^*j_n(\del)}(u_{(n-1)\del}),
\ \text{provided\ \ $u_{(n-1)\del}\in H_{c^*j_n(\del)}$ and $j_n(\del)>0$.}
\]
\end{lemma}

\begin{proof}
We consider only $n=2$; the proof for $n>2$ is similar.
Fix $\del$ and a solution $(u,\beta)$. Arguing by contradiction,
assume $\rho_{2,\del}>b=b(2,\del,u_\del)=\La_{c^*j_2(\del)}(u_\del)$
(the latter is well defined and finite
by the assumption $u_\del\in H_{c^*j_2(\del)}$).
Hence $\beta((b,\iy)\times[\del,2\del])>0$.
Because $\beta$ does not charge $\R\times\{\del\}$, it follows that
$\theta:=\beta((b,\iy)\times(\del,2\del])>0$.

Let us show using \eqref{018}(ii) that there exists $t\in(\del,2\del]$
such that $U(b,t)=0$. If this statement is false, namely $U(b,t)>0$ for all
$t\in(\del,2\del]$, then
\[
\beta(U>0)\ge\int_{(b,\iy)\times(\del,2\del]}1_{\{U(x,t)>0\}}\beta(dx,dt)
\ge\int_{(b,\iy)\times(\del,2\del]}1_{\{U(b,t)>0\}}\beta(dx,dt)=\theta>0,
\]
contradicting \eqref{018}(ii).

Fix such $t\in(\del,2\del]$.
We now appeal to identity \eqref{012} with $\tau=\del$.
Denoting the last term there by
$h(y)=\int_{\R\times[\del,t]}\frs_{t-s}(x,y)\beta(dx,ds)$,
\[
\int_{-\iy}^bh(y)dy\le\int_{-\iy}^\iy h(y)dy
=\int_{[\del,t]} e^{\kap(t-s)}dJ_s
\le\int_{[\del,2\del]} e^{\kap(2\del-s)}dJ_s=j_2(\del).
\]
Since $U(b,t)=0$, we have for the first term in \eqref{012},
\[
\int_{-\iy}^b\int_\R \frs_{t-\del}(x,y)u(x,\del)dx\,dy\le\int_{-\iy}^b h(y)dy\le j_2(\del).
\]
Using \eqref{92}, $\int_{-\iy}^x \frs_{t-\del}(x,y)dy\ge\int_{-\iy}^x \frp_{t-\del}(x,y)dy\ge c_*$ for all $x$, hence
\begin{align*}
j_2(\del)&\ge\int_{x\in(-\iy,b]}\int_{y\in(-\iy,b]}\frs_{t-\del}(x,y)u(x,\del)dy\,dx\\
&\ge c_*\int_{-\iy}^b u(x,\del)dx=c_*c^*j_2(\del)=2j_2(\del),
\end{align*}
a contradiction due to the assumption $j_2(\del)>0$.
This proves the claim.
\end{proof}

\begin{proof}[Proof of Proposition \ref{prop2}]
The proof of \eqref{062} is similar to the proof of
the analogous statement from Proposition \ref{prop1}.
Also as in that proof, the $L_1$ norm of the solution and
the barrier at $n\del$ are equal.
As for \eqref{r3}, fix $0<t_0<T$. Denoting
$\chi_n=1_{(n-1)\del<t_0}+e^{-\Del^5/\del}$,
it follows from \eqref{016} that
\[
\ell_{n,\DEL}\le e^{\kap\del}\ell_{n-1,\DEL}+j_n(\del)\chi_n.
\]
By induction,
$\ell_{n,\DEL}\le e^{n\kap\del}\sum_{i=1}^n j_i(\del)\chi_i$.
Using $j_n(\del)\le e^{\kap\del}(J_{n\del}-J_{(n-1)\del})$
and $n\del\le T$,
\[
\ell_{n,\DEL}\le e^{\kap(T+\del)}(J_{t_0+\del}+e^{-\Del^5/\del}J_T),
\]
as claimed.

We turn to the main assertion.
Assume that $\Del<\eps_0/2$.
Arguing by induction, assume that $f\LE g \MOD \ell_{n-1,\DEL}$, where
\[
f=u^{(\DEL,-)}_{(n-1)\del},
\qquad
g=u_{(n-1)\del},
\]
when $n-1\ge1$ and $f=g=\xi_0$ when $n-1=0$.
Write $C$, $P$ and $S$ for $C^{(\DEL)}_n$, $P^{(\del)}_n$ and $S_\del$,
respectively.
Then $u^{(\DEL,-)}_{n\del}=CPSf$,
and by Lemma \ref{lem1}, $u_{n\del}=PSg-h$, where
\[
h(y)=\int_{\R\times[(n-1)\del,n\del]}\frs_{n\del-s}(x,y)\beta(dx,ds).
\]
By Lemma \ref{lem3}(i), $PSf\LE PSg$ mod $e^{\kap\del}\ell_{n-1,\DEL}$.
If $(n-1)\del<t_0$ we therefore have, for any $r$,
\begin{align*}
u^{(\DEL,-)}_{n\del}[r,\iy)=CPSf[r,\iy)&\le PSf[r,\iy)\\
&\le PSg[r,\iy)+e^{\kap\del}\ell_{n-1,\DEL}
\\
&\le PSg[r,\iy)-h[r,\iy)+\|h\|_1+e^{\kap\del}\ell_{n-1,\DEL}
\\
&=u_{n\del}[r,\iy)+j_n(\del)+e^{\kap\del}\ell_{n-1,\DEL}\\
&=u_{n\del}[r,\iy)+\ell_{n,\DEL},
\end{align*}
which gives the claimed estimate.

In what follows, $(n-1)\del\ge t_0$. In particular, $n\ge2$.
In view of the lower bound on $m_t=\|u(\cdot,t)\|_1$, $t\in[0,T]$
and the continuity of $J$, we may assume that $\del$ is so small that
the condition $u_{(n-1)\del}\in H_{c^*j_n(\del)}$ holds for all
$n\ge 2$, $n\del\le T$. As a result, the bound
asserted in Lemma \ref{lem2} is valid provided merely
that $j_n(\del)>0$.
Moreover, by Lemma \ref{lem1} there exists a constant $c_\iy$ such that
for any solution $(u,\beta)$, $\|u(\cdot,t)\|_\iy<c_\iy$, $t\in[t_0,T]$.
Using the induction assumption and Lemma \ref{lem3}(iv),
\[
u^{(\DEL,-)}_{n\del}=CPSf\LE CPSg \MOD e^{\kap\del}\ell_{n-1,\DEL}.
\]
Denote $w=PSg$. Suppose
\begin{equation}\label{013}
Cw\LE w-h \MOD \eps,
\text{ where $\eps=\ell_{n,\DEL}-e^{\kap\del}\ell_{n-1,\DEL}
=e^{-\Del^5/\del}j_n(\del)$}.
\end{equation}
Then $u^{(\DEL,-)}_{n\del}\LE w-h=u_{n\del}\MOD \ell_{n,\DEL}$,
which completes the proof.
It remains to show \eqref{013}.

First, if $j_n(\del)=0$ then $\eps=0$ and \eqref{013}
holds because $C=\id$, $h=0$, and so $Cw=w=w-h$.
Next assume $j_n(\del)>0$.
Denote $j=j_n(\del)$ and $b=\La_{c^*j_n(\del)}(g)=\La_{c^*j}(g)$.
By Lemma \ref{lem2}, $\rho_{n,\del}\le b$.
Write $h=h_1+h_2$, where
\[
h_1(y)=h(y) 1_{\{y>b+\Del^2\}},
\qquad
h_2(y)=h(y) 1_{\{y\le b+\Del^2\}}.
\]
Because $\rho_{n,\del}\le b$, we have
\begin{align*}
h(y)
=\int_{(-\iy,b]\times[(n-1)\del,n\del]}\frs_{n\del-s}(x,y)\beta(dx,ds).
\end{align*}
W.l.o.g., $e^{\kap\del}<2$, thus $\frs_t\le2\frp_t$ for
$t\le\del$.
Hence in view of Lemma \ref{lem8},
$\frs_t(0,[a,\iy))\le c_3e^{-c_4a^2/t}$ for $a>t^{1/2}>0$,
where $c_3,c_4>0$ depend only on $\hat c_1,\hat c_2$ of the lemma.
This gives
\[
\|h_1\|_1\le c_3\int_{(-\iy,b]\times[(n-1)\del,n\del]}e^{-c_4\Del^4/(n\del-s)}
\beta(dx,ds)
\le c_3e^{-c_4\Del^4/\del}j,
\]
provided $\del<\Del^4$. If we further require $\Del<c_4$
then for all sufficiently small $\del$,
\begin{equation}
\label{053}
\|h_1\|_1\le e^{-\Del^5/\del}j=\eps.
\end{equation}

Recall that $\|h\|_1=j$ and let $q\in(0,1]$ be defined by
$\|h_2\|_1=qj$. Then $q\ge1-e^{-\Del^5\del}$.
Let us argue that it suffices to show
\begin{equation}\label{020}
b+\Del^2=\La_{c^*j}(g)+\Del^2\le \La_\Del(w)
\end{equation}
in order to prove \eqref{013}.
By definition, $h_2$ is supported to the left of $b+\Del^2$.
On the other hand, $C_{\Del,qj}w=w-\tilde h$, where $\|\tilde h\|_1=qj$
and $\tilde h$ is supported to the right of $\La_\Del(w)$.
Thus using $\|\tilde h\|_1=\|h_2\|_1$, it follows from \eqref{020}
that $C_{\Del,qj}w\LE w-h_2= w-h+h_1$. In view of \eqref{053} this gives
\[
C_{\Del,qj}w\LE w-h \MOD \eps.
\]
Because $Cw=C_{\Del,j}w\le C_{\Del,qj}w$ pointwise,
one has $Cw\LE C_{\Del,qj}w$. Hence \eqref{013} follows.

It remains to show \eqref{020}. Because $t_0\le (n-1)\del\le T$,
the bound $\|g\|_\iy\le c_\iy$ is valid.
By making $\Del$ smaller if needed, assume
$2c_\iy\Del^2<\Del/6$. Then for all $\del$ so small that $c^*j<\Del/6$ (simultaneously over $n$),
\begin{equation}
\label{052}
\La_{\Del/3}(g)\ge \La_{c^*j}(g)+2\Del^2.
\end{equation}
Next we argue that for all small $\del$,
\begin{equation}\label{021}
\La_{2\Del/3}(Sg)\ge \theta:=\La_{\Del/3}(g)-\Del^2.
\end{equation}
To show this we must to show $(Sg)(-\iy,\theta]\le2\Del/3$.
Let $g=g_1+g_2=\vhat C_{\Del/3}g+C_{\Del/3}g$.
Because $\|g_1\|_1=\Del/3$,
we have $\|Sg_1\|_1\le e^{\kap\del}\Del/3\le\Del/2$, and
\[
(Sg)(-\iy,\theta]=(Sg_1)(-\iy,\theta]+(Sg_2)(-\iy,\theta]
\le \Del/2+\|g_2\|_1e^{\kap\del}\frp_\del(x_0,(-\iy,\theta]),
\]
where $x_0=\La_{\Del/3}(g)$, owing to the fact that
$\frp_\del(y,(-\iy,\theta])$ is monotone decreasing in $y$ for $y>\theta$.
Recalling that $\|g\|_1=\|u_{(n-1)\del}\|_1=m_{(n-1)\del}\le\|m\|^*_T$
and using again Lemma \ref{lem8},
the last term in the above display is bounded by
\[
2\hat c_1\|m\|^*_T\int_{-\iy}^{x_0-\Del^2}\del^{-1/2}e^{-\hat c_2(x_0-y)^2\del^{-1}}dy,
\]
which is smaller than
$\Del/6$ for all sufficiently small $\del$. This shows
$(Sg)(-\iy,\theta]\le2\Del/3$ hence \eqref{021}.

For $\pi:=S*\al(\cdot,n\del;(n-1)\del)$
we have $\|\pi\|_1\le e^{\kap\del}(I_{n\del}-I_{(n-1)\del})$.
Hence for all small $\del$, $\|\pi\|_1<\Del/3$. As a result,
\[
\La_\Del(w)=\La_\Del(PSg)=\La_\Del(Sg+\pi)
\ge \La_{2\Del/3}(Sg).
\]
Combining this with \eqref{052} and \eqref{021} gives \eqref{020},
and the proof is complete.
\end{proof}

\subsection{Proof of uniqueness}\label{sec33}

The last step is showing that the lower and upper barriers
become close upon taking $\del\to0$ then $\Del\to0$ and finally $t_0\to0$.

\begin{proposition}\label{prop3}
Fix $0<t_0<T$. Let $\Del_0$ and $\del_0=\del_0(\Del_0)$ be as in
Proposition \ref{prop2}. Then
for $\Del\in(0,\Del_0)$, $\del\in(0,\del_0)$ and $n\in\N$, $n\del\le T$, one has
\[
u^{(\del,+)}_{n\del}\LE u^{(\DEL,-)}_{n\del} \MOD e^{n\kap\del}\Del.
\]
\end{proposition}

\begin{proof}
By induction. Recall $u^{(\pm)}_0=\xi_0$.
Assume
$u^{(\del,+)}_{(n-1)\del}\LE u^{(\DEL,-)}_{(n-1)\del} \MOD e^{(n-1)\kap\del}\Del$. Then
\[
P^{(\del)}_nS_\del u^{(\del,+)}_{(n-1)\del}\LE P^{(\del)}_nS_\del u^{(\DEL,-)}_{(n-1)\del} \MOD e^{n\kap\del}\Del,
\]
where, for $n-1=0$, this is true because both sides of the inequality are equal,
and otherwise this is a consequence of the induction assumption and Lemma \ref{lem3}(i),
recalling that the $L_1$ norm of the upper and lower barriers are equal
for each $n$. For the same reason, Lemma \ref{lem3}(vi) also applies, and gives
\begin{align*}
C^{(\del)}_nP^{(\del)}_nS_\del u^{(\del,+)}_{(n-1)\del}
\LE C^{(\DEL)}_nP^{(\del)}_nS_\del u^{(\DEL,-)}_{(n-1)\del} \MOD e^{n\kap\del}\Del,
\end{align*}
that is,
$u^{(\del,+)}_{n\del}\LE u^{(\DEL,-)}_{n\del} \MOD e^{n\kap\del}\Del$.
This completes the proof.
\end{proof}

\begin{proof}[Proof of Theorem \ref{th1}]
Once uniqueness is established for the $u$ component
of the solution $(u,\beta)$, uniqueness of the $\beta$ component
follows from \eqref{054}.
By Remark \ref{rem2} it suffices to prove uniqueness of solutions
$(u,\beta)$ in which $u$ is the version given by Lemma \ref{lem1}.
To show uniqueness of the $u$ component,
argue by contradiction and assume that $(u^i,\beta^i)$, $i=1,2$ are two
solutions where $u^1$ and $u^2$ are distinct.
Then there exist $t>0$ and $r\in\R$ such that,
say, $u^1_t[r,\iy)<u^2_t[r,\iy)$. Fix such $t$ and $r$.
Denote $\del_n=tn^{-1}$ for $n\in\N$. Let $0<t_0<t$.
Then by Propositions \ref{prop1} and \ref{prop2}, for every small $\Del>0$
there exists $n_0$ such that for $n>n_0$,
\[
u^{(\Del,\del_n,-)}_t[r,\iy)-\ell_{n,\Del,\del_n}\le u^1_t[r,\iy)< u^2_t[r,\iy)\le u^{(\Del,+)}_t[r,\iy).
\]
By Proposition \ref{prop3},
\[
u^{(\Del,+)}_t[r,\iy)\le u^{(\Del,\del_n,-)}_t[r,\iy)+e^{\kap t}\Del.
\]
Using these two inequalities and then the bound from \eqref{r3},
\[
0<u^2_t[r,\iy)-u^1_t[r,\iy)\le e^{\kap t}\Del+\ell_{n,\Del,\del_n}
\le e^{\kap t}\Del+e^{\kap(t+\del_n)}(J_{t_0+\del_n}+e^{-\Del^5/\del_n}J_t).
\]
On taking $n\to\iy$, then $\Del\downarrow0$ and finally
$t_0\downarrow0$, the expression on the right converges to zero, a contradiction.
\end{proof}

\section{Convergence}\label{sec4}
\beginsec

In this section the convergence result is proved based on the FBP uniqueness,
yielding the proof of Theorem \ref{th2}.
Throughout, the assumptions of Theorem \ref{th2} hold.
The main steps are as follows.
In Lemma \ref{lem4}, the normalized processes are shown
to satisfy a version of equation \eqref{r4} with an error term.
Lemma \ref{lem6} establishes tightness of these processes.
Existence of a measurable density for any limit point of
the sequence $\bar\xi^N$ is shown in Lemma \ref{lem5}.
In Lemma \ref{lem7}, the final, crucial step shows
that the complementarity condition is preserved under the limit.

\subsection{Limit laws and the parabolic equation}
\label{sec41}

This subsection contains the first three of the aforementioned
steps toward convergence.
Some notation used here is as follows.
Let $J^{\#,N}_t=\beta^N(\R\times[0,t])$ denote the counting
process for removals, and note that, by construction,
$J^{\#,N}_t=J^N_t$ holds
on the event $\{\inf_{s\le t}m^N_s\ge1\}$ (recall the remark after \eqref{r14}).
Let
\[
Y^N_t=\sum_{i\in\calS}1_{\{\sig^i\le t\}}(\pi^i_{t\w\tau^i}-\pi^i_{\sig^i})
\]
be the number of births during $[0,t]$, and let its macroscopic counterpart
be given by $Y_t=\kap\int_0^tm_sds$.
Denote $R^N_T=N+I^N_T$ and note that
the set $\calR^N_T=\{(j,0):j\le R^N_T\}$ consists of all
root particles appearing by time $T$.

Recall that a solution to \eqref{018}(i, iv) is defined via \eqref{054}.
The relation of the particle system to this equation is established
by showing that if $(\xi,\beta)$ is a limit point of $(\bar\xi^N,\bar\beta^N)$ then
\begin{equation}\label{054+}
-\int_0^\iy\lan\pl_t\ph+\calL\ph,\xi_t\ran dt
=\int_{\R}\ph(\cdot,0)d\xi_0+\int_{\R\times\R_+}\ph d\al
-\int_{\R\times\R_+}\ph d\beta.
\end{equation}

\begin{lemma}\label{lem4}
i.
Let $\ph\in C_c^\iy(\R\times\R_+)$ and let $T$ be such that $\ph(\cdot,t)=0$ for all $t\ge T$. Then
\[
-\int_0^\iy\lan(\pl_t\ph+\calL\ph)(\cdot,t),\bar\xi^N_t\ran dt
=\lan\ph(\cdot,0),\bar\xi^N_0\ran
+\int_{\R\times\R_+}\ph d\bar\al^N
-\int_{\R\times\R_+}\ph d\bar\beta^N+\bar M^N_T,
\]
where $\bar M^N_t$ is an $\{\calF^N_t\}$-martingale starting at zero,
with quadratic variation
\begin{equation}\label{019}
[\bar M^N]_t
\le cN^{-1}\Big(\int_0^t\bar m^N_sds+\bar Y^N_t\Big),
\end{equation}
where $c$ depends only on $\ph$ and $\frc$.

ii. One has
\[
\bar m^N_t=1+\kap\int_0^t\bar m^N_sds
+\bar I^N_s-\bar J^{\#,N}_t+\bar M^{\#,N}_t,
\]
where, with $c$ as above, $\bar M^{\#,N}_t$ is an $\{\calF^N_t\}$-martingale
starting at zero, and
\[
[\bar M^{\#,N}]_t\le c N^{-1}\bar Y^N_t.
\]

iii. As $N\to\iy$, $(\bar m^N,\bar Y^N,\bar I^N,\bar J^N,\bar J^{\#,N})
\to(m,Y,I,J,J)$ in probability.

iv. Suppose that $(\xi_0,\xi,\al,\beta)$ is a subsequential limit
of $(\bar\xi^N_0,\bar\xi^N,\bar\al^N,\bar\beta^N)$. Then the former
tuple satisfies \eqref{054+} a.s.
\end{lemma}

\begin{proof}
i.
Note that $\sig^i$ and $\tau^i$ are $\{\calF^N_t\}$-stopping times
and recall that $B^i$ and $\hat\pi^i$ are martingales on this
filtration.
By It\^o's lemma, for each $i\in\calS$, on the event $\{t\ge\sig^i\}$,
\begin{align}\label{31}
\ph(X^i_{t\w\tau^i},t\w\tau^i)&=
\ph(x^i,\sig^i)
+\int_{\sig^i}^{t\w\tau^i}(\pl_t\ph+\frb\pl_x\ph+\fra\pl_x^2\ph)(X^i_s,s)ds
+\int_{\sig^i}^{t\w\tau^i}(\frc\pl_x\ph)(X^i_s,s)dB^i_s
\notag
\\
&=
\ph(x^i,\sig^i)+\int_{\sig^i}^{t\w\tau^i}(\pl_t\ph+\frb\pl_x\ph+\fra\pl_x^2\ph)(X^i_s,s)ds
+M^{N,i,1}_t,
\end{align}
where
\[
M^{N,i,1}_t
=1_{\{\sig^i\ge t\}}\int_{\sig^i}^{t\w\tau^i}(\frc \pl_x\ph)(X^i_s,s)dB^i_s.
\]
Given $i$, the sum of
evaluations of $\ph$ over birth location-epochs of particles born directly
from particle $i$ between time $0$ and $t$ is given by
\[
1_{\{t\ge\sig^i\}}\int_{\sig^i}^{t\w\tau^i} \ph (X^i_s,s)d\pi^i_s.
\]
Summing this expression over $i\in\calS$ gives
the sum of evaluations of $\ph$ over all birth location-epochs
during that time interval, i.e.,
\begin{align*}
\sum_{i=(j,k)\in\calS^N_t,\, k\ge1}\ph(x^i,\sig^i)
&= \sum_{i\in\calS}1_{\{t\ge\sig^i\}}\int_{\sig^i}^{t\w\tau^i} \ph (X^i_s,s)d\pi^i_s\\
&= \sum_{i\in\calS}1_{\{t\ge\sig_i\}}\Big[\int_{\sig^i}^{t\w\tau^i} \ph (X^i_s,s)\kap ds
+M^{N,i,2}_t\Big],
\end{align*}
where the index set on the left corresponds to births by time $t$, and
\[
M^{N,i,2}_t=1_{\{t\ge\sig^i\}}\int_{\sig^i}^{t\w\tau^i} \ph (X^i_s,s)
d\hat\pi^i_s.
\]
Therefore, summing \eqref{31} over all $i$ such that $\sig^i\le t$
and normalizing gives
\begin{align*}
\int_{\R\times[0,t]}\ph d\bar\beta^N
&=\int_\R\ph(\cdot,0)d\bar\xi^N_0
+\int_{\R\times[0,t]}\ph d\bar\al^N
+\int_0^t\int_\R(\pl_t\ph+\calL\ph)(x,t)\bar\xi^N_t(dx)dt
+\bar M^N_t,
\end{align*}
where
\[
\bar M^N_t=N^{-1}M^N_t,
\qquad
M^N_t=\sum_{i\in\calS}(M^{N,i,1}_t+M^{N,i,2}_t).
\]
Take $t=T$ and replace the
integration range $\R\times[0,T]$ to $\R\times\R_+$ recalling
that $\ph(x,t)=0$ for $t>T$.
The bound \eqref{019}
follows from $[M^{N,i,1}]_t\le 1_{\{t\ge\sig_i\}}\|\frc\|^2_\iy\|\pl_x\ph\|^2_\iy(t\w\tau^i-\sig^i)$
and $[M^{N,i,2}]_t\le 1_{\{t\ge\sig_i\}}\|\ph\|^2_\iy(\pi^i_{t\w\tau^i}-\pi^i_{\sig^i})$
and the identities
\[
\sum_{i\in\calS}1_{\{t\ge\sig^i\}}(t\w\tau^i-\sig^i)
=\int_0^t\xi^N_s(\R)ds=\int_0^tm^N_sds
\]
and
\begin{equation}\label{r26}
\sum_{i\in\calS}1_{\{t\ge\sig^i\}}(\pi^i_{t\w\tau^i}-\pi^i_{\sig^i})=Y^N_t.
\end{equation}

ii.
We have $m^N_t=N+I^N_t+Y^N_t-J^{\#,N}_t$
by the definition of these processes. By \eqref{r26},
\begin{equation}\label{r27}
Y^N_t=\kap\int_0^tm^N_sds+\kap M^{\#,N}_t,
\qquad\text{where}
\qquad
M^{\#,N}_t=\sum_{i\in\calS}1_{\{t\ge\sig^i\}}\int_{\sig^i}^{t\w\tau^i}d\hat\pi^i_s.
\end{equation}
The quadratic variation bound follows as in (i).

iii.
Fix $T$. Recall $\eps_0$ from Assumption \ref{assn1}.
Consider the $\{\calF^N_t\}$-stopping time
\[
\theta^N=\inf\{t:\bar I^N_t\ge I_T+1 \text{ or } \bar m^N_t\le \eps_0/2\}.
\]
By (ii) and Gronwall's lemma,
$\E[\bar m^N_{t\w\theta^N}] \le c$, where $c=c(T)$.
Hence by \eqref{r27}, $\E[\bar Y^N_{t\w\theta^N}]\le c$.
Going back to (ii) gives that $\E[[\bar M^{\#,N}]_{T\w\theta^N}]\to0$,
hence $\|\bar M^{\#,N}\|^*_{T\w\theta^N}\to0$ in probability.
For $t\le \theta^N$, one has $J^{\#,N}_t=J^N_t$.
Thus using \eqref{055} and again Gronwall's lemma one has
$\|\bar m^N-m\|^*_{T\w\theta^N}\to0$ in probability.
By the definition of $\eps_0$, this shows that
$\PP(\theta^N<T)\to0$. Because $T$ is arbitrary, this proves
that $\bar m^N\to m$ and $\bar J^{\#,N}\to J$ in probability.
As a result, by \eqref{r27}, $\bar Y^N\to Y$ in probability.

iv. In view of (i), it suffices to show that $\|\bar M^N\|^*_T\to0$ in probability.
However, this is an immediate consequence of \eqref{019} and (iii).
\end{proof}

The following point is used in the proof of the next two lemmas.
One can construct an additional particle system in which there are no removals,
based on the same stochastic primitives as in the original system,
except $\tilde\eta^l=\iy$ for all $l$ in place of the original removal attempt times $\eta^l$.
In this particle system we use tilde notation for all the model ingredients,
as in $\tilde x^i$, $\tilde X^i$, $\tilde\sig^i$,
with one exception: instead of $\tilde\xi^N$ we write $\zeta^N$
(and $\bar\zeta^N$ for its normalized version).
Thus $\tilde J^N=0$, $\tilde\beta^N=0$ and  $\tilde\tau^i=\iy$ for all $i$.
The tilde system dominates the original system in several ways.
For example, it is easy to see by a simple coupling
that for every $A\in\calB(\R)$ and $t\ge0$,
the random variable $\xi^N_t(A)$ is stochastically dominated by $\zeta^N_t(A)$.

\begin{lemma}\label{lem6}
The sequence of laws of $(\bar\xi^N,\bar\beta^N)$, $N\in\N$,
is tight. For every subsequential limit $(\xi,\beta)$,
one has $\PP(\xi\in C(\R_+,\calM_+(\R)))=1$.
\end{lemma}

\begin{proof}
Both the $J_1$ topology over $D(\R_+,\R)$ and
the topology of local weak convergence we gave the space
$\calM_{+,\rm loc}(\R\times\R_+)$ are defined by convergence over
finite time intervals.
Hence in this proof we fix $T$ and consider all processes
(respectively, measures) defined on $\R_+$ (on subsets of $\R\times\R_+$)
as if they are defined on $[0,T]$ (on subsets of $\R\times[0,T]$),
and with a slight abuse of notation still use the same notation.
For example, $\beta^N$ will denote the restriction of the original
random measure $\beta^N$ to subsets of $\R\times[0,T]$.

Moreover, recalling that $\bar\al^N(\R\times[0,T])=\bar I^N_T\to I_T$
in probability, we may and will assume w.l.o.g.\ that
the injections are truncated when their number reaches $c_IN$,
where $c_I=I_T+1$. Hence, for all $N$ and $t\in[0,T]$, $I^N_t\le c_IN$ a.s.
Similarly, the removal measure is assumed w.l.o.g.\ to be truncated
when $J^N$ reaches $c_JN$, $c_J=J_T+1$.

Tightness may be argued separately for each component. Starting with $\bar\beta^N$,
for $r>0$, denote
$\mathbb{B}_{r,T}=[-r,r]\times[0,T]$, and for $n\in\N$, set
\[
K_n(r)=\{\gam\in\calM_+(\R\times[0,T]): |\gam|\le c_J,\,
\gam(\mathbb{B}_{r,T}^c)<n^{-1}\}.
\]
Then by Prohorov's theorem (see
\cite[Theorem~A2.4.I]{daley03} for a version for finite measures),
for any $n_0$ and sequence $\{r_n\}$, the closure of
$K_{\ge n_0}(\{r_n\}):=\cap_{n\ge n_0}K_n(r_n)$ as a subset of
$\calM(\R\times[0,T])$ endowed with the topology of weak convergence,
is compact.
Suppose we show that for every $n$ there exists $r_n$ such that
\begin{equation}\label{r13}
\liminf_N\PP(\bar\beta^N\in K_n(r_n))\ge1-2^{-n}.
\end{equation}
Given $\eps>0$ let $n_0$ be such that
$\sum_{n\ge n_0}2^{-n}<\eps$. Then
\[
\liminf_N\PP(\bar\beta^N\in \oo{K_{\ge n_0}(\{r_n\})})>1-\eps.
\]
This would show that $\bar\beta^N$ are tight.

Toward showing \eqref{r13}, note that the removal space-time location of a particle
is a point on the graph of the potential trajectory of that particle, hence
\begin{equation}\label{r10}
\beta^N(\mathbb{B}_{r,T}^c)
\le U^N_r:=\sum_{i\in\calS^N_T}1_{\{\|X^i\|^*_{[\sig^i,T]}>r\}}.
\end{equation}
If we let
\[
\tilde U^N_r=\sum_{i\le\tilde\calS^N_T}1_{\{\|\tilde X\|^*_{[\tilde\sig^i,T]}>r\}},
\]
then for each $N$ and $r$, the random variable $U^N_r$
is stochastically dominated by $\tilde U^N_r$,
as follows by a simple coupling between the two particle
systems.
Recalling that $|\bar\beta^N|\le c_J$, this shows
\begin{equation}\label{r15}
\PP(\bar\beta^N\in K_n(r)^c)
=\PP(\bar\beta^N(\mathbb{B}_{r,T}^c)\ge n^{-1})
\le\PP(\tilde U^N_r\ge Nn^{-1}).
\end{equation}

In the tilde system, the collection of family members descending from
a root particle $i=(j,0)$ up to time $T$
is denoted by $\tilde\calS^{N,i}_T$, in accordance with \eqref{r20}. Let
\begin{align*}
\tilde\calF^N_*
&=\sig\{I^N_T,\,(\tilde x^i,\tilde\sig^i):i\in\calR^N_T\}
\\
&=\sig\{I^N_T,\,(x^i,\sig^i):i\in\calR^N_T\},
\end{align*}
where the equality follows by construction.
For $i\in\calR^N_T$ and $\hat i\in\tilde\calS^{N,i}_T$, let $\tilde X^{i,\hat i}_t$,
$t\in[\tilde\sig^i,T]$ denote the trajectory formed by $\tilde X^{\hat i}_t$
during its lifetime, and by the trajectories of its ancestors
prior to its birth time (here as well $\tilde\sig^i$ can be replaced by
$\sig^i$).
Recalling that the motion and branching mechanisms
are independent of the initial configuration and injection measure,
using the many-to-one lemma \cite{harris2017}, we have
\begin{equation}\label{r16}
\E\Big[\sum_{\hat i\in\tilde\calS^{N,i}_T}1_{\{\|\tilde X^{i,\hat i}\|^*_{[\sig^i,T]}>r\}}
\Big|\tilde\calF^N_*\Big]
=e^{\kap(T-\sig^i)}\E[1_{\{\|\tilde X^i\|^*_{[\sig^i,T]}>r\}}|\tilde\calF^N_*],
\qquad i\in\calR^N_T.
\end{equation}
Let $X_t$ solve \eqref{407} for $t\in[0,T]$, and denote the stochastic
integral term by $C_t=\int_0^t\frc(X_\theta)dB_\theta$.
Then $\lan C\ran=\int_0^\cdot\frc(X_s)^2ds$, and
by time change for continuous martingales,
\begin{equation}\label{r12}
C_t=\check B_{\lan C\ran_t}
\end{equation}
where $\check B_s=C_{\tau(s)}$ is a BM, and
$\tau(s)=\inf\{t\ge0: \lan C_t\ran> s\}$. 
By the boundedness of the coefficients $\frb,\frc$, this gives
$\|X-x\|^*_T\le c_1(1+\|\check B\|^*_{c_1})$,
where $c_1$ depends only on $T$ and the coefficients.
Let $\check B^i$ denote the BM, constructed via time change as above,
corresponding to $\tilde X^i$, $i\in\calR^N_T$, and note that for such
$i$, $\tilde x^i=x^i$. Then we have shown
$\|\tilde X^i-x^i\|^*_{[\sig^i,T]}\le c_1(1+\|\check B^i\|^*_{c_1})$.
(Note that the BM $\check B^i$ are not in general mutually independent).
Hence given $n$, there exists $r'=r'_n$ such that
\begin{equation}\label{r17}
\PP(\|\tilde X^i-x^i\|^*_{[\sig^i,T]}>\frac{r'}{2}|\tilde\calF^N_0)
\le\PP(c_1(1+\|\check B^1\|^*_{c_1})>\frac{r'}{2})\le (2n(1+c_I)e^{\kap T})^{-1}2^{-n},
\quad
i\in\calR^N_T.
\end{equation}

Next, by our assumptions, the normalized configuration
measure of $\{x^i:i\in\calR^N_T\}$, given by
$\bar\xi^N_0+\int_{[0,T]}\bar\al^N(\cdot,dt)$, converges in probability
to a deterministic finite measure on $\R$.
Hence for every $n$ there exists $r''=r''_n$ such that
\[
\lim_N\PP(\#\{i\in\calR^N_T:|x^i|>\frac{r''}{2}\}<N(2e^{\kap T})^{-1}2^{-n}) =1.
\]
Hence recalling $R^N_T\le (1+c_I)N$,
\begin{equation}\label{r18}
\limsup_N\E\Big[\frac{\#\{i\in\calR^N_T:|x^i|>r''/2\}}{N}\Big]
\le (2ne^{\kap T})^{-1}2^{-n}.
\end{equation}

For $n\in\N$ let $r_n=r'_n\vee r''_n$. Then by \eqref{r15},
\eqref{r16}, \eqref{r17} and finally \eqref{r18},
\begin{align}\label{r19}
\notag
\limsup_N\PP(\bar\beta^N\in K_n(r_n)^c)
&\le\limsup_N\frac{n}{N}\E[\tilde U^N_{r_n}]
\\
\notag
&\le\limsup_N\frac{n}{N}\E\sum_{i\in\calR^N_T}\sum_{\hat i\in\tilde\calS^{N,i}_T}
1_{\{\|\tilde X^{i,\hat i}\|^*_{[\sig^i,T]}>r_n\}}
\\
\notag
&\le \limsup_N\frac{ne^{\kap T}}{N}\E\sum_{i\in\calR^N_T}
(1_{\{\|\tilde X^i-x^i\|^*_{[\sig^i,T]}>r_n/2\}}+1_{\{|x^i|>r_n/2\}})
\\
&\le
2^{-n}.
\end{align}
This shows \eqref{r13} hence the tightness of the laws of $\bar\beta^N$.

Denote by $d_{\rm L}$ the Levy-Prohorov metric on $\calM_+(\R)$,
which is compatible
with weak convergence on $\calM_+(\R)$.
We will use the notation $w_T(\cdot,\cdot)$ for both $(\R,|\cdot|)$
and $(\calM_+(\R),d_{\rm L})$.
The argument for $\bar\xi^N$ is based on showing
(i) for every $\eps>0$ there exists a compact set $K\subset\calM_+(\R)$
such that $\liminf_N\inf_{t\in[0,T]}\PP(\bar\xi^N_t\in K)>1-\eps$;
and (ii)
for every $\eps>0$ there exists $\del>0$ such that
\[
\limsup_N \PP(w_T(\bar\xi^N,\del)>\eps)<\eps.
\]
Once these two properties are proved it will follow that
$\bar\xi^N$ is a relatively compact sequence
\cite[Corollary 3.7.4 (p.\ 129)]{ethkur}.
Because we use $w$ rather than $w'$
\cite[(3.6.2) (p.\ 122)]{ethkur}, this will in fact establish $C$-tightness,
proving the second statement.

To show (i), let $c_2=\|m\|^*_T+1$. By Lemma \ref{lem4}(iii),
$\|\bar m^N\|^*_T=\sup_{t\in[0,T]}|\bar\xi^N_t|\le c_2$ w.h.p.
Denoting
\[
\hat K_n(r)=\{\gam\in\calM_+(\R):|\gam|\le c_2,\, \gam([-r,r]^c)<n^{-1}\},
\]
it suffices to prove that for every $n$ there exists $r_n$ such that
\[
\liminf_N\inf_t\PP(\bar\xi^N_t\in\hat K_n(r_n))\ge1-2^{-n},
\]
for the same reason given above for \eqref{r13} to be sufficient for
tightness of $\bar\beta^N$. Moreover,
similarly to the estimate \eqref{r10} for $\beta^N$, we have
$\xi^N_t([-r,r]^c)\le U^N_r$ for all $t\in[0,T]$.
Hence the chain of inequalities \eqref{r19} provides a bound also on
$\limsup_N\sup_t\PP(\bar\xi^N_t\in\hat K_n(r_n)^c)$, and (i) follows.

It remains to show (ii).
Let $A^N_t=\{i\in\calS:\sig^i\le t<\tau^i\}$. This is the index set for living particles
at time $t$. For $\del>0$, let $a_\del=\{(s,t)\in[0,T]:0<t-s\le\del\}$.
For $(s,\del)\in a_\del$, the number of particles removed during $(s,t]$
is $J^{\#,N}_t-J^{\#,N}_{s}$.
The number of new particles during this interval is given by
$I^N_t-I^N_{s}+Y^N_t-Y^N_s$.
Hence, denoting symmetric difference by $\triangle$,
\[
\#(A^N_s\triangle A^N_t)\le J^{\#,N}_t-J^{\#,N}_s+I^N_t-I^N_s
+Y^N_t-Y^N_s.
\]
The convergence of $\bar I^N$, $\bar J^{\#,N}$ and $\bar Y^N$
to continuous paths
shows that given $\eps>0$ there exists $\del>0$ such that,
w.h.p., for all $(s,t)\in a_\del$ one has
$\#(A^N_s\triangle A^N_t)\le\eps N/2$.

Next, going back to \eqref{r12} and the notation $\check B^i$,
there exists a constant $c_3\in(0,\iy)$ such that
$w_{[\sig^i,T]}(X^i,\del)\le c_3\del+w_{[\sig^i,c_3T]}(\check B^i,c_3\del)$.
Thus if $p(\eps,\del)=\PP(c_3\del+w_{c_3T}(B,c_3\del)\ge\eps)$
for $B$ a BM, then for $i\in\calS$,
\[
\PP(w_{[\sig^i,T]}(X^i,\del)\ge\eps|\calF^N_{\sig^i})\le p(\eps,\del)
\qquad \text{on } \{\sig^i\le T\}.
\]
Hence
\begin{align*}
\E[\#\{i\in\calS^N_T:w_{[\sig^i,T]}(X^i,\del)\ge\eps\}]
&=\sum_{i\in\calS}\E[ 1_{\{\sig^i\le T\}}\E[ 1_{\{w_{[\sig^i,T]}(X^i,\del)\ge\eps\}} | \calF^N_{\sig^i}]]
\\
&\le \E[\#\calS^N_T]p(\eps,\del)
\\
&\le cNp(\eps,\del),
\end{align*}
where in the last inequality we used the fact that
the expected number of descendants each root particle has by time $T$
is bounded, which along with the truncation convention of $I^N$ gives
$\E[\#\calS^N_T]\le cN$. This gives
\[
\PP(\#\{i\in\calS^N_T:w_{[\sig^i,T]}(X^i,\del)\ge\eps\}>\eps N/2)
\le\frac{cNp(\eps,\del)}{\eps N/2}<\frac{\eps}{4},
\]
where we used the fact that $p(\eps,0+)=0$
and chose $\del=\del(\eps)$ sufficiently small.

Given a set $C\subset\R$ let $C^\eps$ denote its $\eps$-neighborhood.
We have shown the following.
For all $N$ so large that
\[
\PP(\text{for some }
(s,t)\in a_\del,\ A^N_s\triangle A^N_t>\eps N/2)<\frac{\eps}{4},
\]
with probability greater than $1-\eps/2$, for all $(s,t)\in a_\del$,
except for at most $\eps N/2$ particles (removed between $s$ and $t$),
and at most $\eps N/2$ particles (whose displacement exceeds $\eps$),
each particle $i\in A^N_s$ exists in the configuration at time $t$
and travels less than $\eps$ between $s$ and $t$.
Hence, with probability greater than $1-\eps/2$, for any Borel set $C$,
\[
\xi^N_s(C) \le \xi^N_t(C^\eps) + \eps N,
\qquad (s,t)\in a_\del.
\]
Similarly, with probability greater than $1-\eps/2$,
\[
\xi^N_t(C) \le \xi^N_s(C^\eps)+ \eps N,
\qquad (s,t)\in a_\del.
\]
Hence with probability $\ge1-\eps$,
\[
d_{\rm L}(\bar\xi^N_s,\bar\xi^N_t)\le\eps,
\qquad (s,t)\in a_\del.
\]
This shows that
\[
\limsup_N \PP(w_T(\bar\xi^N,\del)>\eps)\le\eps,
\]
and the proof is complete.
\end{proof}

Recall that
\[
v(y,t)=S_t\xi_0(y)+S*\al(y,t)=\int_{\R}\frs_t(x,y)\xi_0(dx)
+\int_{\R\times[0,t]}\frs_{t-s}(x,y)\al(dx,ds),
\]
and set $\zeta_t(dx)=v(x,t)dx$, $t>0$, and $\zeta_0(dx)=\xi_0(dx)$.

\begin{lemma}\label{lem5}
Let $(\xi,\beta)$ be a subsequential limit of $(\bar\xi^N,\bar\beta^N)$.
Then there exists an event $\Om_1\in\calF$ of full measure and
a $\calB(\R)\otimes\calB(\R_+)\otimes\calF$-measurable
function $u(x,t,\om)$ such that for every $(t,\om)\in(0,\iy)\times\Om_1$,
$u(\cdot,t,\om)$ is a density of $\xi_t(\cdot,\om)$ with respect to
the Lebesgue measure on $\R$.
Moreover, for almost all $\om$ one has
\[
u(x,t,\om)\le v(x,t) \qquad\text{for all } (x,t)\in\R\times(0,\iy).
\]
\end{lemma}

\begin{proof}
Let $\calQ$ denote the set of bounded open intervals
$(a,b)\subset\R$. Fix $T>0$.
In the first step we show that for $t\in(0,T]$, $Q\in\calQ$ and $\gam>0$,
\begin{equation}\label{050+}
\PP(\bar\zeta^N_t(Q)>\zeta_t(Q)+\gam)\to0.
\end{equation}
For a first moment calculation, using the many to one lemma as before,
\begin{align*}
\EE[\bar\zeta^N_t(Q)|\tilde\calF^N_*]
&=N^{-1}\E[\sum_{i\in\calR^N_t}\sum_{\hat i\in\calS^{N,i}_t}1_{\{\tilde X^{i,\hat i}_t\in Q\}}|\tilde\calF^N_*]
\\
&=N^{-1}\sum_{i\in\calR^N_t}
e^{\kap(t-\sig^i)}\PP(\tilde X^i_t\in Q|\tilde\calF^N_*)
\\
&=\bar \Th^N:=N^{-1}\sum_{i\in\calR^N_t}\theta^i(N),
\end{align*}
where $\theta^i(N)=\frs_{t-\sig^i}(x^i,Q)$. Now, for $\eps\in(0,t)$,
\begin{align}\label{r22}\notag
\bar\Th^N
&=\int_{\R}\frs_t(x,Q)\bar\xi^N_0(dx)+\int_{\R\times[0,t]}\frs_{t-s}(x,Q)\bar\al^N(dx,ds)
\\ \notag
&\le
V^N(t,\eps,Q):=
\int_{\R}\frs_t(x,Q)\bar\xi^N_0(dx)+\int_{\R\times[0,t-\eps]}\frs_{t-s}(x,Q)\bar\al^N(dx,ds)
\\
&
\hspace{8em}
+e^{\kap T}(\bar I^N_t-\bar I^N_{t-\eps}).
\end{align}
On $\R\times[0,t-\eps]$, $(x,s)\mapsto\frp_{t-s}(x,Q)$
is bounded and continuous \cite[Theorem 1.2.1]{stroock}.
Since $\al$ does not charge $\R\times\{t-\eps\}$,
$V^N(t,\eps,Q)$ converges in probability to
\[
V(t,\eps,Q):=\int_{\R}\frs_t(x,Q)\xi_0(dx)
+\int_{\R\times[0,t-\eps]}\frs_{t-s}(x,Q)\al(dx,ds)+c(I_t-I_{t-\eps}).
\]
We may adopt the truncation convention from the proof of Lemma \ref{lem6}.
Thus $\bar I^N_T$ and $|\bar\al^N|$ are bounded,
and we have by bounded convergence $\E[V^N(t,\eps,Q)]\to V(t,\eps,Q)$,
which, upon taking $\eps\to0$, gives
\[
\limsup\E[\bar\Th^N]\le
\int_QS_t\xi_0(y)dy+\int_QS*\al(y,t)dy = \zeta_t(Q).
\]
Moreover, dropping the last term in \eqref{r22} gives
a lower bound on $\bar\Th^N$, hence similarly
$\liminf\E[\bar\Th^N]\ge\zeta_t(Q)$, giving
\begin{equation}\label{r23}
\lim\E[\bar\zeta^N_t(Q)]=\lim\E[\bar\Th^N]=\zeta_t(Q).
\end{equation}
Similarly,
\begin{equation}\label{r21}
\limsup\E[(\bar\Th^N)^2]\le\zeta_t(Q)^2.
\end{equation}
By Lemma \ref{lem1}, $\zeta_t(Q)<\iy$.
For a second moment calculation, if $i\in\calR^N_t$ then
\[
\E\Big[\sum_{\hat i_1,\hat i_2\in\tilde\calS^{N,i}_t}1_{\{\tilde X^{\hat i_1}_t\in Q,
\tilde X^{\hat i_2}_t\in Q\}}\Big|\tilde\calF^N_*\Big]\le \E[(\tilde Z^{N,i}_T)^2|\tilde\calF^N_*]\le c,
\]
whereas if $i_1,i_2\in\calR^N_T$ are distinct, using conditional independence
and the many-to-one lemma,
\[
\E\Big[\sum_{\hat i_1\in\tilde\calS^{N,i_1}_T, \hat i_2\in\tilde\calS^{N,i_2}_T}
1_{\{\tilde X^{\hat i_1}_t\in Q,\tilde X^{\hat i_2}_t\in Q\}}\Big|\tilde\calF^N_*\Big]
=e^{\kap(t-\sig^{i_1}+t-\sig^{i_2})}\PP(\tilde X^{i_1}_t\in Q|\tilde\calF^N_*)
\PP(\tilde X^{i_2}_t\in Q|\tilde\calF^N_*).
\]
Therefore
\begin{align*}
\E[\bar\zeta^N_t(Q)^2|\tilde\calF^N_*]
&\le cN^{-1}+(\bar\Th^N)^2,
\end{align*}
and using now \eqref{r21}, $\limsup\E[(\bar\zeta^N_t(Q))^2]\le\zeta_t(Q)^2$.
In view of \eqref{r23}, this gives that
$\lim\text{var}(\bar\zeta^N_t(Q))=0$.
By \eqref{r23} this shows \eqref{050+}.

Because $\bar\zeta^N_t$ dominates $\bar\xi^N_t$, \eqref{050+}
holds for the latter as well. In the next
step it is shown that for every $t\in(0,T]$ and $Q\in\calQ$
there is a full-measure event on which
\begin{equation}\label{050}
\xi_t(Q)\le \zeta_t(Q).
\end{equation}
Since along a subsequence one has
$\bar\xi^N\To\xi$ and the latter has continuous sample paths,
one also has $\bar\xi^N_t\To\xi_t$.
Using Skorohod's representation we may assume w.l.o.g.\ that
$\bar\xi^N_t\to\xi_t$ a.s., and since $Q$ is open, we have
$\liminf\bar\xi^N_t(Q)\ge\xi_t(Q)$ a.s. Hence
with $\tilde c=\zeta_t(Q)+\gam$,
\[
\PP(\xi_t(Q)>\tilde c)
\le\PP(\liminf\bar\xi^N_t(Q)>\tilde c)
\le\E\liminf1_{\{\bar\xi^N_t(Q)>\tilde c\}}\le\liminf\PP(\bar\xi^N_t(Q)>\tilde c)
=0,
\]
where the second inequality is by the lower semicontinuity of
$x\mapsto1_{\{x>\tilde c\}}$ and the third is by Fatou's lemma.
This shows $\PP(\xi_t(Q)>\zeta_t(Q)+\gam)=0$ for every $\gam>0$
hence \eqref{050}.

Let $\tilde\calQ\subset\calQ$ be the set of open intervals $(a,b)$
with $a,b\in\Q$. Then there exists an event
$\Om_0$ of full measure
on which for all $t\in(0,T]\cap\Q$ and all $Q\in\tilde\calQ$,
\eqref{050} holds.
Using the continuity of $t\mapsto\xi_t$ we have that, on $\Om_0$,
\eqref{050} holds for all $(t,Q)\in(0,T]\times\tilde\calQ$. The last
assertion can be extended to $(0,T]\times\calQ$ by taking
$\tilde\calQ\ni Q_n\uparrow Q\in\calQ$.
It follows that on an event of full measure, for all $t\in(0,T]$,
\begin{equation}\label{r25}
\xi_t(A)\le\zeta_t(A),\qquad A\in\calB(\R),
\end{equation}
\cite[Corollary 2, p.\ 169]{bil2}; in particular, $\xi_t(dx)\ll dx$.
Since $T$ is arbitrary, this statement holds with $(0,T]$ replaced by $(0,\iy)$.
Assuming $\xi=0$ outside
the full measure event, we finally obtain that
for every $(t,\om)\in(0,\iy)\times\Om$, $\xi_t(dx,\om)\ll dx$.
We can now appeal to \cite[Theorem 58 in Chapter V (p.\ 52)]{del-mey-82}
and the remark that follows.
The measurable spaces denoted in \cite{del-mey-82}
by $(\Om,\calF)$ and $(T,\calT)$ are taken to be
$(\R,\calB(\R))$ and $((0,\iy)\times\Om,\calB((0,\iy))\otimes\calF)$, respectively.
According to this result there exists
a $\calB(\R)\otimes\calB((0,\iy))\otimes\calF$-measurable function
$u(x,t,\om)$, such that for every $(t,\om)\in(0,\iy)\times\Om$,
$u(\cdot,t,\om)$ is a density of $\xi_t(dx,\om)$ with respect to $dx$.

For the last assertion of the lemma, note that by \eqref{r25},
modifying $u$ into $u\w v$ still gives a density of $\xi_t(\cdot,\om)$.
\end{proof}

\subsection{Limit laws and the complementarity condition}\label{sec42}

Here we prove that the complementarity condition carries over to the limit.

\begin{lemma}\label{lem7}
i.
Let $u\in L_{1,\rm loc}(\R_+,L_1)\cap L_{\iy,\loc}((0,\iy),L_\iy)$ be positive.
Denote $U(x,t)=\int_{-\iy}^xu(y,t)dy$. Let
$\beta\in\calM_{+,\rm loc}(\R\times\R_+)$, $\beta(\R\times\{0\})=0$.
Then the following two conditions are equivalent:
\begin{equation}
\label{c1}
\beta(U>0)=0
\end{equation}
\begin{equation}
\label{c2}
\calI_r:=\int_{[r,\iy)\times\R_+}U(r,t)\beta(dx,dt)=0
\qquad \text{for all } r\in\R.
\end{equation}
ii.
Let $(\xi,\beta)$ be a subsequential limit of
$(\bar\xi^N,\bar\beta^N)$. Then a.s.,
\begin{equation}\label{025}
\calI_r(\xi,\beta):=\int_{[r,\iy)\times\R_+}\xi_t(-\iy,r]\beta(dx,dt)=0\qquad \text{for all } r\in\R.
\end{equation}
iii.
Consequently, if $u$ is the density from Lemma \ref{lem5}
(defined arbitrarily on $\R\times\{0\}$)
and $U(x,t,\om)=\int_{-\iy}^xu(y,t,\om)dy$ then
$(U,\beta)$ satisfy \eqref{c1} a.s.
\end{lemma}

\begin{proof}
i.
To show that \eqref{c1} implies \eqref{c2}, write
\[
\calI_r\le\int_{[r,\iy)\times\R_+}U(x,t)\beta(dx,dt)
\le\int_{\R\times\R_+}U(x,t)\beta(dx,dt)=0.
\]
For the converse,
assume \eqref{c1} is false. Then there exists $\del>0$, $\beta(U>\del)>0$.
Since $\beta$ does not charge $\R\times\{0\}$,
there exist $0<t_1<t_2<\iy$ and finite $a<b$
such that, with $K=[a,b]\times[t_1,t_2]$,
\[
\beta(K\cap \{U>\del\})>0.
\]
With $c$ an upper bound on $u$ in $K$, and $\eps=\frac{\del}{2c}\w(b-a)$,
\[
U(y,t)-U(x,t)\le \frac{\del}{2} \qquad x,y\in[a,b],\ 0\le y-x\le\eps,\ t\in[t_1,t_2].
\]
Moreover, there exists $r\in[a,b]$ such that $\beta(L)>0$ where
\[
L=[r,r+\eps]\times[t_1,t_2]\cap\{U>\del\}.
\]
Hence for $(x,t)\in L$,
\[
U(r,t)\ge U(x,t)-\frac{\del}{2}\ge\del-\frac{\del}{2}=\frac{\del}{2}.
\]
Thus
\[
\calI_r\ge\int_L U(r,t)\beta(dx,dt)\ge\frac{\del}{2}\beta(L)>0,
\]
showing that \eqref{c2} is false.

ii.
Let $\Sig$ be the collection of tuples
$\sig=(r,t_1,t_2,\eta,\del)\in\Q^5$,
$0<t_1<t_2$, $\eta>0$, $\del>0$. We show that for every
$\sig\in\Sig$, $\PP(\Om_\sig)=0$ where
\[
\Om_\sig=\{\inf_{t\in[t_1,t_2]}\xi_t(-\iy,r]>\eta,\, \beta((r,\iy)\times(t_1,t_2))>\del\}.
\]
Fix $\sig=(r,t_1,t_2,\eta,\del)$. If $\PP(\Om_\sig)>0$ then by
the weak convergence $(\bar\xi^N,\bar\beta^N)\To(\xi,\beta)$,
the a.s.\ continuity of the limit $\xi$, and the fact that,
for every $t>0$, the measure $\xi_t$ has no atoms,
one must have for all large $N$,
\[
\PP\Big(\inf_{t\in[t_1,t_2]}\bar\xi^N_t(-\iy,r]>\eta/2,
\,\bar\beta^N(r,\iy)\times(t_1,t_2)>\del/2\Big)>0.
\]
However, by the construction of the particle system,
for every $r$, a removal never occurs at a location $>r$
at a time when there are particles at location $\le r$.
Hence the above probability is zero for all $N$.
This shows $\PP(\Om_\sig)=0$.
Consequently, $\PP(\cup_\Sig\Om_\sig)=0$.

Next consider the event
\[
\Om^0=\{\text{there exists $r\in\R$ such that $\calI_r(\xi,\beta)>0$}\}.
\]
On this event there exists $r\in\R$ and $0<s_1<s_2<\iy$
such that
\[
\calI^*:=\int_{[r,\iy)\times(s_1,s_2)}\xi_t(-\iy,r]\beta(dx,dt)>0.
\]
Consider $\Q\ni r_n\uparrow r$. Recalling that the density
$u(\cdot,\cdot,\om)$ is bounded by $v$ and denoting
$\gam_1=\sup_{(x,t)\in[r-1,r]\times(s_1,s_2)}v(x,t)$,
$\gam_2=\beta(\R\times(s_1,s_2))$,
\begin{align*}
\calI_{r_n}(\xi,\beta)
&\ge\int_{[r_n,\iy)\times(s_1,s_2)}\xi_t(-\iy,r_n]\beta(dx,dt)
\\
&\ge\int_{[r,\iy)\times(s_1,s_2)}(\xi_t(-\iy,r]-\xi_t(r_n,r])\beta(dx,dt)
\\
&\ge\calI^*-\sup_{t\in(s_1,s_2)}\xi_t(r_n,r]\gam_2
\\
&\ge
\calI^*-(r-r_n)\gam_1\gam_2>0
\end{align*}
for large $n$. This shows that
on $\Om^0$ there exists $r\in\Q$ such that $\calI_r(\xi,\beta)>0$.

Next, the condition $\calI_r(\xi,\beta)>0$ (with $r\in\Q$) implies that
there exists $\eta\in\Q\cap(0,1)$ such that
$\int_{A_\eta}a_tdb_t>0$ where we denote $a_t=\xi_t(-\iy,r]=\xi_t(-\iy,r)$,
$b_t=\beta(r,\iy)\times[0,t]$ and $A_\eta=\{t:a_t>2\eta\}$.
The trajectory $t\mapsto a_t$ is continuous on $(0,\iy)$
(using the fact that $\xi\in C(\R_+,\calM_+(\R))$ and
that for each $t$, $\xi_t$ has no atoms).
Hence there exists an interval $(t_1,t_2)\subset A_\eta$,
with $t_1,t_2\in\Q$, such that $\int_{(t_1,t_2)}a_tdb_t>0$.
Consequently,
\[
\xi_t(-\iy,r]=a_t\ge2\eta>\eta
\]
on $[t_1,t_2]$, while
\[
\beta((r,\iy)\times(t_1,t_2))=\int_{(t_1,t_2)}db_t>0.
\]
This shows that $\PP(\Om^0)\le\PP(\cup_{\sig\in\Sig}\Om_\sig)=0$.

iii.
The final assertion follows from the first two as soon as
these conditions are verified:
$\|u(\cdot,t,\om)\|_\iy$ is locally bounded for $t\in(0,\iy)$, and
$\beta(\R\times\{0\})=0$. The former follows from Lemmas \ref{lem5}
and \ref{lem1}(iii), by which $u(\cdot,\cdot,\om)\le v$ and
$v\in L_{\loc,\iy}((0,\iy),L_\iy)$.
The latter follows from Lemma \ref{lem4}(iii) by which
$\bar J^{\#,N}\to J$ and the assumption $J_0=0$.
\end{proof}

\subsection{Proof of Theorem \ref{th2}}\label{sec43}

In view of the tightness stated in Lemma \ref{lem6}
and the uniqueness of solutions to \eqref{018}
stated in Theorem \ref{th1},
it suffices to show that whenever $(\xi_0,\xi,\al,\beta,J)$ is a subsequential limit
of $(\bar\xi^N_0,\bar\xi^N,\bar\al^N,\bar\beta^N,\bar J^N)$,
and $u$ the corresponding density from Lemma \ref{lem5}, one has that,
a.s., $(u,\beta)$ is a solution to \eqref{018}.

That $u\in L_{1,\rm loc}(\R_+,L_q)$ for $q\in(1,\iy)$
follows from Lemma \ref{lem5},
which states that $u\le v$,
and Lemma \ref{lem1}(i), by which
$v\in L_{1,\rm loc}(\R_+,L_q)$ for all $q\in(1,\iy)$.
Since by Lemma \ref{lem4}(iii) $\bar J^{\#,N}\to J$,
we have $\beta(\R\times[0,t])=J_t<\iy$ for all $t$, showing that
$\beta\in\calM_{+,\rm loc}(\R\times\R_+)$.
Thus to show that $u$ is a weak $L_q$-solution to \eqref{018}(i),
it remains to show that \eqref{054} holds. This is indeed the case
by Lemma \ref{lem4}(iv), in view of the relation $\xi_t(dx)=u(x,t)dx$.

Finally, Lemma \ref{lem7} shows that condition \eqref{018}(ii) holds,
and for condition \eqref{018}(iii) we have just provided a proof.
This shows that, a.s., $(u,\beta)$ is a solution to \eqref{018},
and completes the proof.
\qed

{\bf Acknowledgment}
The author would like to thank Kavita Ramanan for valuable discussions and comments.
This research was supported by ISF grant 1035/20.

\bibliographystyle{is-abbrv}

\bibliography{rev1}

\begin{thebibliography}{10}
\ifx \showCODEN  \undefined \def \showCODEN #1{CODEN #1}  \fi
\ifx \showISBN   \undefined \def \showISBN  #1{ISBN #1}   \fi
\ifx \showISSN   \undefined \def \showISSN  #1{ISSN #1}   \fi
\ifx \showLCCN   \undefined \def \showLCCN  #1{LCCN #1}   \fi
\ifx \showPRICE  \undefined \def \showPRICE #1{#1}        \fi
\ifx \showURL    \undefined \def \showURL {URL }          \fi
\ifx \path       \undefined \input path.sty               \fi
\ifx \ifshowURL \undefined
     \newif \ifshowURL
     \showURLtrue
\fi

\bibitem{ama01}
H.~Amann.
\newblock Linear parabolic problems involving measures.
\newblock {\em Real Academia de Ciencias Exactas, Fisicas y Naturales. Revista.
  Serie A, Matematicas}, 95\penalty0 (1):\penalty0 85--119, 2001.

\bibitem{ber21}
J.~Berestycki, {\'E}.~Brunet, J.~Nolen, and S.~Penington.
\newblock A free boundary problem arising from branching {B}rownian motion with
  selection.
\newblock {\em Transactions of the American Mathematical Society}, 374\penalty0
  (09):\penalty0 6269--6329, 2021.

\bibitem{ber22bee}
J.~Berestycki, {\'E}.~Brunet, J.~Nolen, and S.~Penington.
\newblock Brownian bees in the infinite swarm limit.
\newblock {\em The Annals of Probability}, 50\penalty0 (6):\penalty0
  2133--2177, 2022.

\bibitem{ber19}
J.~Berestycki, {\'E}.~Brunet, and S.~Penington.
\newblock Global existence for a free boundary problem of {F}isher--{KPP} type.
\newblock {\em Nonlinearity}, 32\penalty0 (10):\penalty0 3912, 2019.

\bibitem{bil2}
P.~Billingsley.
\newblock {\em Probability and Measure}.
\newblock John Wiley, New York, third edition, 1995.

\bibitem{bru06}
E.~Brunet, B.~Derrida, A.~H. Mueller, and S.~Munier.
\newblock Noisy traveling waves: effect of selection on genealogies.
\newblock {\em Europhysics Letters}, 76\penalty0 (1):\penalty0 1, 2006.

\bibitem{bru07}
{\'E}.~Brunet, B.~Derrida, A.~H. Mueller, and S.~Munier.
\newblock Effect of selection on ancestry: an exactly soluble case and its
  phenomenological generalization.
\newblock {\em Physical Review E}, 76\penalty0 (4):\penalty0 041104, 2007.

\bibitem{dem19}
M.~Cabezas, A.~Dembo, A.~Sarantsev, and V.~Sidoravicius.
\newblock Brownian particles with rank-dependent drifts: Out-of-equilibrium
  behavior.
\newblock {\em Communications on Pure and Applied Mathematics}, 72\penalty0
  (7):\penalty0 1424--1458, 2019.

\bibitem{de-masi-top}
G.~Carinci, A.~De~Masi, C.~Giardin{\`a}, and E.~Presutti.
\newblock Hydrodynamic limit in a particle system with topological
  interactions.
\newblock {\em Arabian Journal of Mathematics}, 3\penalty0 (4):\penalty0
  381--417, 2014.

\bibitem{de-masi-book}
G.~Carinci, A.~De~Masi, C.~Giardin{\`a}, and E.~Presutti.
\newblock {\em Free {B}oundary {P}roblems in {P}{D}{E}s and {P}article
  {S}ystems, {\rm {S}pringer briefs in mathematical physics}}, volume~12.
\newblock Springer, 2016.

\bibitem{daley03}
D.~J. Daley and D.~Vere-Jones.
\newblock {\em An Introduction to the Theory of Point Processes: Volume I:
  Elementary Theory and Methods}.
\newblock Springer, 2003.

\bibitem{de-masi-excl}
A.~De~Masi, P.~A. Ferrari, and E.~Presutti.
\newblock Symmetric simple exclusion process with free boundaries.
\newblock {\em Probability Theory and Related Fields}, 161\penalty0
  (1-2):\penalty0 155--193, 2015.

\bibitem{de-masi-nbbm}
A.~De~Masi, P.~A. Ferrari, E.~Presutti, and N.~Soprano-Loto.
\newblock Hydrodynamics of the {$N$}-{BBM} process.
\newblock In {\em International workshop on Stochastic Dynamics out of
  Equilibrium}, pages 523--549. Springer, 2017.

\bibitem{de-masi-nbbm2}
A.~De~Masi, P.~A. Ferrari, E.~Presutti, and N.~Soprano-Loto.
\newblock Non local branching {B}rownian motions with annihilation and free
  boundary problems.
\newblock {\em Electronic Journal of Probability}, 24, 2019.

\bibitem{delarue22}
F.~Delarue, S.~Nadtochiy, and M.~Shkolnikov.
\newblock Global solutions to the supercooled {S}tefan problem with blow-ups:
  regularity and uniqueness.
\newblock {\em Probability and Mathematical Physics}, 3\penalty0 (1):\penalty0
  171--213, 2022.

\bibitem{del-mey-82}
C.~Dellacherie and P.-A. Meyer.
\newblock {\em Probabilities and {P}otential {B}. Transl. by {JP} {W}ilson}.
\newblock North Holland, Amsterdam, 1982.

\bibitem{dur11}
R.~Durrett and D.~Remenik.
\newblock Brunet--{D}errida particle systems, free boundary problems and
  {W}iener--{H}opf equations.
\newblock {\em The Annals of Probability}, 39\penalty0 (6):\penalty0
  2043--2078, 2011.

\bibitem{ethkur}
S.~Ethier and T.~Kurtz.
\newblock {\em Markov Processes: Characterization and Convergence}.
\newblock Wiley, New York, 1986.

\bibitem{folland95}
G.~B. Folland.
\newblock {\em Introduction to {P}artial {D}ifferential {E}quations}.
\newblock Princeton university press, 1995.

\bibitem{harris2017}
S.~C. Harris and M.~I. Roberts.
\newblock The many-to-few lemma and multiple spines.
\newblock In {\em Annales de l'Institut Henri Poincar{\'e}, Probabilit{\'e}s et
  Statistiques}, volume~53, pages 226--242. Institut Henri Poincar{\'e}, 2017.

\bibitem{klu23}
A.~Klump.
\newblock The inverse first-passage time problem as hydrodynamic limit of a
  particle system.
\newblock {\em Methodology and Computing in Applied Probability}, 25\penalty0
  (1):\penalty0 42, 2023.

\bibitem{landim98}
C.~Landim, S.~Olla, and S.~Volchan.
\newblock Driven tracer particle in one dimensional symmetric simple exclusion.
\newblock {\em Communications in mathematical physics}, 192\penalty0
  (2):\penalty0 287--307, 1998.

\bibitem{landim06}
C.~Landim and G.~Valle.
\newblock A microscopic model for {S}tefan's melting and freezing problem.
\newblock {\em The Annals of Probability}, 34\penalty0 (2):\penalty0 779--803,
  2006.

\bibitem{mai16}
P.~Maillard.
\newblock Speed and fluctuations of $n$-particle branching {B}rownian motion
  with spatial selection.
\newblock {\em Probability Theory and Related Fields}, 166\penalty0
  (3):\penalty0 1061--1173, 2016.

\bibitem{stroock}
D.~W. Stroock.
\newblock {\em Partial Differential Equations for Probabalists}.
\newblock Cambridge University Press, 2008.

\end{thebibliography}

\end{document}